\documentclass[11pt]{amsart}
\usepackage{amsmath}
\usepackage{amssymb}
\usepackage{amsfonts}
\usepackage{amsthm}\usepackage{mathrsfs}
\usepackage{comment}
\usepackage[normalem]{ulem}
\usepackage{enumerate}

\usepackage{geometry}
\usepackage{needspace}
\usepackage{mathtools}
\mathtoolsset{showonlyrefs}
\usepackage{hyperref}
\hypersetup{colorlinks=true,
linktoc=all,linkcolor=black,
citecolor=black}
\usepackage{tikz}
\usepackage{bbm}
\usepackage{mathptmx}

\newcommand{\bea}{\begin{eqnarray}}
\newcommand{\eea}{\end{eqnarray}}

\newcommand{\D}{\mathbb{D}}

\def\textmatrix#1&#2\\#3&#4\\{\bigl({#1 \atop #3}\ {#2 \atop #4}\bigr)}
\def\dispmatrix#1&#2\\#3&#4\\{\left({#1 \atop #3}\ {#2 \atop #4}\right)}
\newcommand{\be}{\begin{equation}}
\newcommand{\ee}{\end{equation}}
\newcommand{\ben}{\begin{eqnarray*}}
\newcommand{\een}{\end{eqnarray*}}

\newcommand{\bi}{\begin{itemize}}
\newcommand{\ei}{\end{itemize}}

\newcommand{\T}{\mathbb{T}}
\newcommand{\imag}{\mathrm{i}}
\newcommand{\e}{\mathrm{e}}

\newcommand{\Nop}{\mathrm{N}}

\newcommand{\re}{\mathrm{Re}}

\newcommand{\C}{\mathbb{C}}

\newcommand{\diff}{\mathrm{d}}

\newcommand{\Ordo}{\mathrm{O}}
\newcommand{\ordo}{\mathrm{o}}

\newcommand{\dA}{\mathrm{dA}}

\newcommand{\Bergman}{\mathrm{P}}

\theoremstyle{definition}

\theoremstyle{plain}

\newtheorem{thm}{Theorem}[subsection]
\newtheorem{cor}[thm]{Corollary}
\newtheorem{lem}[thm]{Lemma}
\newtheorem{prop}[thm]{Proposition}
\theoremstyle{definition}

\newtheorem{rem}[thm]{Remark}

\newtheorem{conj}[thm]{Conjecture}

\numberwithin{equation}{subsection}
\numberwithin{figure}{subsection}

\title{Geometric folding and maximal asymptotic variance}
\author[Haakan Hedenmalm]{Haakan Hedenmalm}
\address{\textbf{Haakan Hedenmalm:}\newline
Department of Mathematics and Computer Science,
St. Petersburg State University, St. Petersburg, Russia\newline
Beijing Institute for Mathematical Sciences and Applications (BIMSA),
Huairou District, Beijing 101408, China\newline
Department of Mathematics and Statistics,
University of Reading, Reading, United Kingdom\newline
Department of Mathematics,
University of Aveiro, Aveiro, Portugal}
\email{haakan00@gmail.com}
\date{22 November 2026}
\usepackage{graphicx}
\usepackage{placeins}


\begin{document}
\begin{abstract}
We develop a finite-scale geometric folding method for the
nonlinearity of a schlicht function, combining a signed
critical-point count with hyperbolic zero-packing rigidity and a
quantitative saddle--maximum pairing. A target cutoff and a
cylindrical root transform allow the resulting strict estimate to
hold uniformly over the unrestricted schlicht class. In particular,
we obtain $\Sigma_{\max}^2<1$ for its maximal asymptotic variance.

We also investigate the regular hyperbolic $\{12,3\}$ zero lattice
as a concrete comparison model. Its optimized cell discrepancy
$0.0523107153\ldots$ is remarkably close to the conjectural lower
benchmark $0.0522117855\ldots$ arising from a circular Neumann cell.
A second-variation calculation gives a strictly positive quadratic
variation for the circular cell along nontrivial smooth
area-preserving perturbations fixing its zero, lending credence
to the possibility that it is optimal among such Neumann cells.
If the circular cell were optimal among all such Neumann cells
of equal area, and if optimal zero configurations admitted
sufficiently well-behaved gradient-flow partitions of the disk,
we would obtain the explicit lower bound
$\rho_{\mathbb H}\ge 0.0522117855\ldots$.
We also discuss Beltrami-generated boundary domains as exploratory
models, without claiming their universal spectral optimality.
Independently of the strict folding estimate, a holomorphic
Beltrami-family argument gives the numerically explicit lower bound
$6(1-\rho_{\mathbb H})\ge5.686\ldots$ for maximal
Schwarzian asymptotic variance.
\end{abstract}
\maketitle

\section{Introduction}
\subsection{Main results}
Let $\mathcal S$ denote the schlicht class of normalized univalent
maps $f:\D\to\C$ with $f(0)=0$ and $f'(0)=1$.
Throughout the paper, for a holomorphic function $h$ on $\D$ and
$0<r<1$, we write $h_r(z):=h(rz)$.  In particular,
$f_r(z)=f(rz)$.  For $f\in\mathcal S$, the normalized dilation
$r^{-1}f_r$ belongs to $\mathcal S$.  Thus
$(f_r)'(z)=r f'(rz)$, while $(r^{-1}f_r)'(z)=f'(rz)$.
We choose the branch of $\log f'$ which vanishes at the origin, and put
\[
L(r):=\log\frac1{1-r^2},\qquad 0<r<1.
\]
We use normalized arclength $\diff s(\zeta)=\diff\theta/(2\pi)$ on
the unit circle $\mathbb T$, where $\zeta=\e^{\imag\theta}$.
For $0<r<1$, we call the circle mean of
$|\log f'(r\zeta)|^2$ the \emph{finite-step variance}.
Our principal result estimates this quantity uniformly in $f$
before taking the limit $r\to1^-$.

\begin{thm}
\label{thm:main-finite-scale}
There exists a universal constant $0<\Sigma_{\mathrm{fold}}^2<1$
such that, as $r\to1^-$,
\begin{equation}\label{eq:strict-main-intro}
 \sup_{f\in\mathcal S}\int_{\mathbb T}
 |\log f'(r\zeta)|^2\,\diff s(\zeta)
 \le \Sigma_{\mathrm{fold}}^2 \,L(r)
       +\mathrm O\!\big(\sqrt{L(r)}\big).
\end{equation}
The implied constant and the lower admissible radius are absolute.
\end{thm}

For clarity, we define the \emph{unrestricted maximal asymptotic variance}
by
\[
\Sigma_{\max}^2
:=\limsup_{r\to1^-}\sup_{f\in\mathcal S}
\frac1{L(r)}\int_{\mathbb T}
|\log f'(r\zeta)|^2\,\diff s(\zeta).
\]
The theorem immediately gives the following strict bound for the
unrestricted maximal asymptotic variance of the schlicht class.

\begin{cor}
\label{cor:cylindrical-maximal}
\[
 \Sigma_{\max}^2\le\Sigma_{\mathrm{fold}}^2<1.
\]
\end{cor}

The finite-radius form of the theorem is important because its error is
uniform over the unrestricted schlicht class. The proof begins with a
corresponding estimate for the bounded subclass
\[
\mathcal S_M:=\{f\in\mathcal S:\,|f(z)|\le M\text{ for }z\in\D\},
\qquad M\ge1,
\]
where $M\ge1$ is necessary for this class to be nonempty.
Writing $\Nop(f)=f''/f'$ for the nonlinearity (or pre-Schwarzian) and
$\dA(z)=\pi^{-1}\diff x\diff y$ for normalized area measure, the
quantitative folding surplus yields, for some absolute constant $\delta>0$,
\[
\int_{|z|<r}(1-|z|^2)|\Nop(f)(z)|^2\,\dA(z)
\le (1-\delta)L(r)+C_\delta(1+\log M),\qquad f\in\mathcal S_M.
\]
The Littlewood--Paley comparison and cylindrical root transform
then remove the image-bound cost without losing the strict leading
coefficient.

\subsection{Background and proof strategy}
Brennan's conjecture and the geometry of harmonic measure provide
historical context for the boundary problem; see Makarov
\cite{MakarovHausdorff} and Carleson--Makarov
\cite{CarlesonMakarovBrennan,CarlesonMakarovAggregation}.
Building on the September 2026 OpenAI preprints
\cite{OpenAIBrennan2026,OpenAIStrict2026}, Binder, Ivrii,
Prause, and Smirnov~\cite{BinderIvriiPrauseSmirnov2026}
developed the signed-area folding method. It replaces the direct
study of a conformal map by the critical-point topology of an
associated nonholomorphic map.

Our proof localizes that signed-area inequality by a radial taper,
uses hyperbolic zero-packing rigidity to obtain a quantitative
saddle--maximum surplus, and then removes the image bound by a
cylindrical root transform. The distinction between uniform and
individual variances is related to~\cite{HedenmalmAdvances2017},
while the zero-packing input comes from~\cite{HedenmalmAJM}.
For complementary integral-means methods see
\cite{HedenmalmShimorinDuke2005,HedenmalmShimorinPAMS2007,
BaranovHedenmalm2008}. We also use the circle-norm comparison
discussed in~\cite{DebnathHedenmalmDirichlet}.

\subsection{Zero packing, Neumann cells, and candidate domains}
We use the regular hyperbolic $\{12,3\}$ zero lattice as an explicit
benchmark, not a proven optimizer. Gradient-flow partitions of the
hyperbolic modulus suggest Neumann cells and the conjectural bound
$\rho_{\mathbb H}\ge\rho_*=0.0522117855\ldots$.
We also study the associated Beltrami family as an exploratory model
for boundary geometry, in the spirit of random conformal snowflakes
\cite{BeliaevSmirnovSnowflakes}. None of these comparisons is needed
for the strict finite-scale theorem.

\section{Geometric folding and the finite-scale estimate}
The proof combines the signed-area folding identity, a radial taper with
target localization, and zero-packing rigidity. A quantitative pairing
of saddles with higher maxima converts absolute Jacobian mass into a
strict signed surplus. Littlewood--Paley comparison and a cylindrical
root transform then give the unrestricted result.

\subsection{The folding Jacobian and signed counting}
\begin{lem}
Let $\Psi:(0,1]\to\mathbb R$ be a real-valued function of class $C^3$
with $\Psi'(s)>0$ for $0<s\le1$. Put $s=1-|z|^2$ and $a(s)=1/\Psi'(s)$, and define
\[G(z)=f(z)+\frac{a(s)}{\bar z}f'(z),\qquad
 V(z)=\Psi(s)-2\log|f(z)-G(z)|.\]
For $z\ne0$, we then have
\[
\e^VJ_G=W(s)\left\{|\Nop(f)|^2+
\frac{2(1-a'(s))}{a(s)}\operatorname{Re}(z\Nop(f))- R(s)\right\},
\quad W(s):=\e^{\Psi(s)},
\]
where
\[
R(s):=\frac1{1-s}+\frac{2a'(s)}{a(s)}+\frac{(1-s)(2a'(s)-1)}{a(s)^2}.
\]
\end{lem}
\begin{proof}
Direct differentiation gives
\[
\partial_z G=f'(1-a'(s)+a(s)\Nop(f)/\bar z),\quad
\bar\partial_z G=-f'(a(s)+(1-s)a'(s))/\bar z^2.
\]
Moreover,
\[
\e^V=\e^{\Psi(s)} |z|^2/(a(s)^2|f'|^2).
\]  The final step is to expand $|\partial_zG|^2-|\bar\partial_zG|^2$.
\end{proof}

The one-pole Morse count gives the signed folding inequality, provided
the target-localized area integral is absolutely convergent.

\begin{lem}
Suppose $f$ is holomorphic and univalent on a neighborhood of
$\overline\D$, with $f'\ne0$ there. Let $\Psi$ be smooth with
$\Psi'>0$, $(1-s)\Psi''-\Psi'<0$, and
$\Psi(s)\to-\infty$ as $s\to0^+$.
For a regular target $\xi\notin f(\partial\D)$ of the critical-point
map $G$, put
\[
u_\xi(z)=\Psi(1-|z|^2)-2\log|f(z)-\xi|.
\]
If $M_\xi(t)$ and $S_\xi(t)$ count its ordinary maxima and saddles
above a regular level $t$, and $C_\xi(t)$ counts superlevel components,
then
\[
M_\xi(t)-S_\xi(t)=C_\xi(t)-P_\xi\ge0,
\]
where $P_\xi\in\{0,1\}$ counts the possible logarithmic pole.
When such a pole $z_\xi=f^{-1}(\xi)$ exists, we extend $u_\xi$
by $u_\xi(z_\xi)=+\infty$ when forming its superlevel sets.
Consequently, for any nonnegative $\eta\in C_c^\infty(\mathbb C)$
for which the following integral is absolutely convergent,
\[
\int_\D\eta(G(z))\e^{V(z)}J_G(z)\,\dA(z)\le0,
\qquad V(z)=u_{G(z)}(z).
\]
\end{lem}

\begin{proof}
Fix a regular target $\xi$. The boundary decay makes each superlevel set
$\Omega_{\xi,t}=\{u_\xi>t\}$ relatively compact after adjoining the
possible logarithmic pole at height $+\infty$. No component has a
hole: on a bounded complementary hole $H$ one would have $u_\xi\le t$
and $u_\xi=t$ on $\partial H$, contradicting the minimum principle for
superharmonic functions. The pole cannot lie in such a hole.
Thus maxima create components as the level is lowered, while saddles
merge two components. The pole, if present, contributes one component
at every finite level. Hence
$M_\xi(t)-S_\xi(t)=C_\xi(t)-P_\xi\ge0$.

At a critical point, with $H(z,\xi)=\partial_z u_\xi(z)$, implicit
differentiation of $H(z,G(z))=0$ gives
\[
 \det D^2u_\xi(z)=-4|H_\xi(z,\xi)|^2J_G(z)
 =-\frac{4|f'(z)|^2}{|f(z)-\xi|^4}J_G(z).
\]
As $\Delta u_\xi<0$, maxima correspond to $J_G<0$ and saddles
to $J_G>0$. Integrating the component count against $\e^t\diff t$
yields
\[
 \sum_{\mathrm{max}}\e^{u_\xi}-\sum_{\mathrm{saddle}}\e^{u_\xi}
 =\int_{-\infty}^{\infty}\e^t(C_\xi(t)-P_\xi)\,\diff t\ge0.
\]
The pole is not an ordinary maximum. Applying the area formula first
on compact exhaustions avoiding the poles of $G$ and the source
boundary, then passing to the limit by the assumed absolute
integrability, proves the claim. Sard's theorem excludes the null
set of nonregular targets.
\end{proof}

\subsection{The radial cutoff and target localization}
Write $s(r)=1-r^2$ and $L(r)=\log\frac1{s(r)}$.
Fix a smooth concave function
$Q:\mathbb R\to\mathbb R$ satisfying
\[
Q(t)=0\quad(t\le0),\qquad Q'(t)=-1\quad(t\ge1),
\qquad -1\le Q'(t)\le0,\quad Q''(t)\le0.
\]
Such a function is obtained by integrating twice a nonnegative smooth
function supported in $(0,1)$ with integral one. This choice is
independent of $r$, $M$, and $f$.

Recall that $s=1-|z|^2$. The resulting smooth radial weight is
\[
W_r(s):=s\exp\!\left(Q\!\left(\log\frac{s(r)}{s}\right)\right),\qquad 0<s\le1.
\]
For the differential calculations, put
\[
L:=\log\frac1s,\qquad Q_r(L):=Q(L-L(r)),\qquad
B_r(L):=1-\frac{\diff Q_r}{\diff L}(L),\qquad
A_r(s):=\frac{s}{B_r(L)}.
\]

\begin{figure}[ht]
\centering
\includegraphics[width=.79\textwidth]{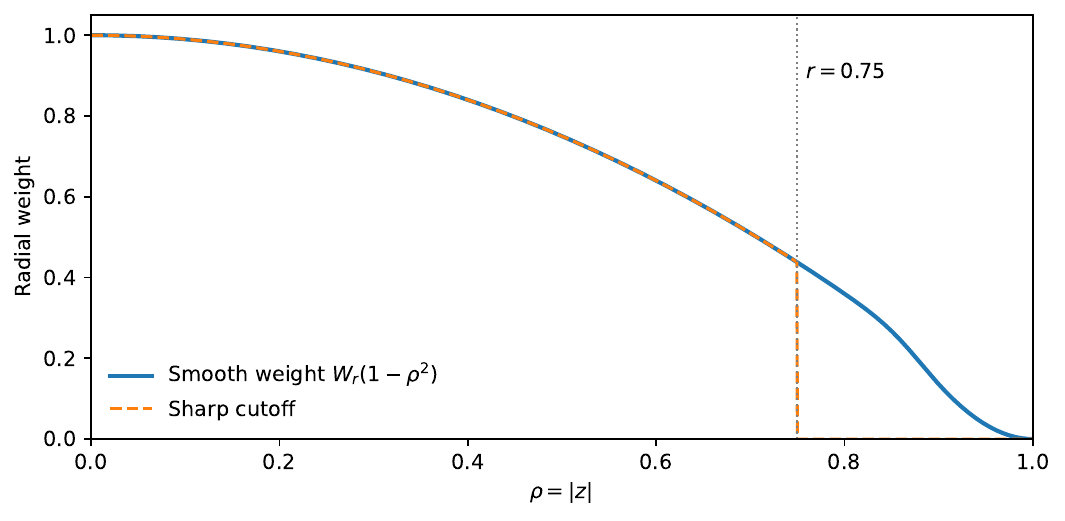}
\caption{The smooth radial weight $W_r(1-\rho^2)$ for one
admissible choice of $Q$, and the sharp cutoff $(1-\rho^2)
\mathbf 1_{\{\rho<r\}}$, for $r=0.75$.
The weights agree for $\rho\le r$, while the smooth profile
tapers toward the unit circle.}
\label{fig:radial-profile}
\end{figure}

Then $\Psi_r(s):=\log W_r(s)$ satisfies $\Psi_r'(s)>0$ and
\[
(1-s)\Psi_r''(s)-\Psi_r'(s)
=-\frac{B_r(L)+(1-s)\dfrac{\diff B_r}{\diff L}(L)}{s^2}<0.
\]
Here $B_r$ is differentiated with respect to $L$, whereas $A_r$
is differentiated with respect to $s$. We define the
\emph{radial Jacobian remainder} by
\[
R_r(s):=\frac1{1-s}+\frac{2A_r'(s)}{A_r(s)}
+\frac{(1-s)(2A_r'(s)-1)}{A_r(s)^2}.
\]
Equivalently, since
\[
A_r'(s)=\frac1{B_r(L)}+
\frac1{B_r(L)^2}\frac{\diff B_r}{\diff L}(L),
\]
we obtain, suppressing the argument $L$ of $B_r$ for readability,
\[
\begin{aligned}
R_r(s)={}&\frac1{1-s}
+\frac2s\left(1+\frac1{B_r}\frac{\diff B_r}{\diff L}\right)\\
&+\frac{1-s}{s^2}\left(2B_r-B_r^2
+2\frac{\diff B_r}{\diff L}\right).
\end{aligned}
\]
We need to compare the radial integrals uniformly.

\begin{lem}
For every fixed $s_0\in(0,1)$, uniformly as $r\to1^-$,
\[
\int_0^{s_0}W_r(s) R_r(s)\,\mathrm ds=L(r)+\mathrm{O}_{s_0}(1).
\]
\end{lem}
\begin{proof}
After the substitution $s=\e^{-L}$, the exact integrand is
\[
\e^{Q_r}\left[(1-s)\left(2B_r-B_r^2+2\frac{\diff B_r}{\diff L}\right)
+2s\left(1+\frac1{B_r}\frac{\diff B_r}{\diff L}\right)
+\frac{s^2}{1-s}\right] \,\diff L,
\]
where $B_r$ and its derivative are evaluated at $L$. For
$L\ge L_0=\log(1/s_0)$, the terms carrying $s$ are uniformly integrable.
Since $\diff Q_r/\diff L=1-B_r$, integration by parts gives
\[
\begin{aligned}
&\int_{L_0}^\infty \e^{Q_r}
\left(2B_r-B_r^2+2\frac{\diff B_r}{\diff L}\right)\,\diff L\\
&\qquad=\int_{L_0}^\infty \e^{Q_r}B_r^2\,\diff L
+2[\e^{Q_r}B_r]_{L_0}^{\infty}.
\end{aligned}
\]
For $L(r)>L_0$, the boundary term equals $-2$. On $[L_0,L(r)]$
we have $Q_r=0$ and $B_r=1$, so the integral contributes
$L(r)-L_0$. The fixed-width transition and exponentially decaying
tail contribute $\mathrm{O}(1)$.
\end{proof}

The target cutoff equals one on the image of the outer annulus and
vanishes near the central pole of $G_r$. Its pullback need not be radial:
in the central region the mixed Jacobian coefficient vanishes, and on
the outer annulus the cutoff is one on complete circles.

The following construction keeps its central cost logarithmic in $M$.

\begin{lem}
Fix once and for all $b_0\in(0,1)$, independently of $M$ and $r$. For $M\ge1$ and $f\in\mathcal S_M$
extending univalently past $\partial\D$, there is a target
cutoff $\eta\in C_c^\infty(\mathbb C)$, with $0\le\eta\le1$,
such that $\eta(G_r)=1$ on $|z|\ge b_0$ and
$\eta(G_r)=0$ on $|z|\le a_M$, for all $r$ sufficiently
close to $1$. The radius $a_M$ can be chosen with
$a_M\ge c_{b_0}/M$, and the central radial contribution satisfies
\[
\int_{|z|<b_0}\eta(G_r(z))W_r(s)R_r(s)\,\dA(z)
\le 2\log M+\mathrm O_{b_0}(1).
\]
\end{lem}
\begin{proof}
The Cauchy estimate and $A_r(s)\le s$ show that
$|G_r(z)|\le M(1+2/b_0)$ for $|z|\ge b_0$.
Choose a smooth cutoff $\eta=1$ on the disk of this radius and $\eta=0$
outside a disk of twice that radius. Its decay begins only beyond
the disk containing the outer-annulus image. The pullback transition
therefore lies on the central sheet. On $|z|\le b_0$,
the folding map is unmodified when $r>b_0$.
For an explicit universal bound, write
$f(z)=z+\sum_{n\ge2}a_nz^n$, with $|a_n|\le n$.
On $|z|\le a_0:=1/4$ the untapered expression satisfies
\[
 G_r(z)-\frac1{\bar z}
 =f(z)-zf'(z)+\frac{f'(z)-1}{\bar z},
\]
so
\[
 \left|G_r(z)-\frac1{\bar z}\right|
 \le\sum_{n\ge2}n(n-1)4^{-n}
     +\sum_{n\ge2}n^2 4^{-(n-2)}
 =\frac{220}{27}<9.
\]
Using the bound $220/27<9$, choose
\[
 K_{b_0}:=1+2/b_0,\qquad
 a_M:=\min\left\{\frac14,\frac1{2K_{b_0}M+10}\right\}.
\]
Since $M\ge1$, we have $a_M\ge c_{b_0}/M$ with
$c_{b_0}=\min\{1/4,1/(2K_{b_0}+10)\}>0$.
Then $|G_r(z)|>2K_{b_0}M$ for $0<|z|\le a_M$,
so $\eta(G_r)=0$ there.
Here $W_rR_r=1/(|z|^2(1-|z|^2))$, so the integral is at most
\[
\int_{a_M<|z|<b_0}
\frac{\dA(z)}{|z|^2(1-|z|^2)}
=\log\frac{b_0^2(1-a_M^2)}
{a_M^2(1-b_0^2)}
\le 2\log M+\mathrm O_{b_0}(1).
\]
\end{proof}

For later use we verify absolute integrability of the signed-area
density with its target cutoff, including the dependence on the
image bound.  Assume that $f$ extends univalently across
$\overline\D$, put $\varepsilon=1-r^2$, and abbreviate
\[
 \mathcal J_r(z):=\e^{u_{G_r(z),r}(z)}J_{G_r}(z).
\]
In the fixed central disk $|z|\le b_0$, the radial profile is
untapered, and the exact Jacobian identity and
$s|\Nop(f)|\le6$ give
\[
 |\mathcal J_r(z)|
 \le \frac1{s|z|^2}+\frac{36}{s},
 \qquad s=1-|z|^2.
\]
Since $\eta(G_r)=0$ on $|z|\le a_M$ and
$a_M\ge c_{b_0}/M$, it follows that
\begin{equation}\label{eq:absolute-central-localization}
 \int_{|z|<b_0}\eta(G_r(z))|\mathcal J_r(z)|\,\dA(z)
 \le C_{b_0}(1+\log M).
\end{equation}
On $|z|\ge b_0$, the quantities $B_r$, $B_r^{-1}$,
$\diff B_r/\diff L$, and $A_r'(s)$ are uniformly bounded,
while $A_r(s)\asymp s$.  The radial Jacobian formula therefore
gives, uniformly in the normalized schlicht map,
\[
 |\mathcal J_r(z)|\le C_{b_0,Q}\frac{W_r(s)}{s^2},
 \qquad |z|\ge b_0.
\]
In the normalized area convention, radial integration is
$\dA=-\diff s$ after angular averaging.  Since $W_r=s$
for $s\ge\varepsilon$ and
$W_r(s)=\e^{1+Q(1)}s^2/\varepsilon$ for
$0<s\le\varepsilon/\e$, splitting the remaining interval
$\varepsilon/\e<s<\varepsilon$ yields
\begin{equation}\label{eq:absolute-outer-localization}
 \int_{|z|\ge b_0}\frac{W_r(s)}{s^2}\,\dA(z)
 \le L(r)+C_{b_0,Q}.
\end{equation}
Consequently
\[
 \int_\D \eta(G_r(z))\,\e^{u_{G_r(z),r}(z)}
             |J_{G_r}(z)|\,\dA(z)
 \le C_{b_0,Q}\bigl(1+L(r)+\log M\bigr)<\infty.
\]
This proves the absolute-integrability hypothesis needed to
pass from compact exhaustions to complete fibers in the signed
area formula.  The bound is independent of the particular map
and is uniform for $r$ sufficiently close to $1$ at the stated
logarithmic scale.

The mixed term is handled without differentiating the target
cutoff.  Its coefficient vanishes identically on $|z|\le b_0$,
where $A_r(s)=s$ and $A_r'(s)=1$.  On $|z|\ge b_0$ the cutoff
is $1$, and the circular mean of
$\operatorname{Re}(z\Nop(f)(z))$ is zero because
$z\Nop(f)$ is holomorphic and vanishes at the origin.
The preceding absolute bounds justify this angular integration.
Thus the mixed term has integral zero and the radial comparison
and central cutoff bounds apply without an additional radiality
assumption on $\eta\circ G_r$.

\subsection{Zero-packing rigidity and the quantitative folding surplus}
\label{subsec:zero-packing-rigidity}

There is a direct connection between the absolute folding Jacobian and
hyperbolic zero packing as studied in~\cite{HedenmalmAJM}. This
connection supplies a quantitative rigidity estimate without a
compactness argument. The strict signed folding inequality follows after combining
this estimate with the complete-fiber Morse pairing below.

Fix $b_0\in(0,1)$, and work on the untapered annulus
$\mathcal A_r=\{z:b_0<|z|<r\}$, where $r>b_0$. Write
\[
 s=1-|z|^2,\qquad g_f(z)=z\Nop(f)(z),\qquad
 j_f(z)=s^2|g_f(z)|^2-1.
\]
In this region the potential is $s/|f(z)-\xi|^2$, its
critical-target map is $G(z)=f(z)+sf'(z)/\bar z$, and the
weighted Jacobian formula reduces to
\[
 \e^{u_{G(z)}(z)}J_G(z)
   =\frac{j_f(z)}{s|z|^2}.
\]
Set
\[
\begin{aligned}
 I_r&:=\int_{\mathcal A_r}\frac{\dA(z)}{s|z|^2},\\
 P_r&:=\int_{\mathcal A_r}\frac{(j_f(z))_+}{s|z|^2}\,\dA(z),\\
 N_r&:=\int_{\mathcal A_r}\frac{(j_f(z))_-}{s|z|^2}\,\dA(z).
\end{aligned}
\]
Here $P_r$ measures the positive, or saddle, mass and $N_r$
measures the negative, or maximum, mass. Notice that
$I_r=L(r)+\mathrm O_{b_0}(1)$.

Recall the hyperbolic discrepancy function
\[
 \Phi_g(z):=\bigl(s|g(z)|-1\bigr)^2
\]
and the minimal discrepancy density $\rho_{\mathbb H}>0$ from
\cite{HedenmalmAJM}. It is defined using the lower limit as
$r\to1^-$ of the infimum of
\[
 \frac1{L(r)}\int_{|z|<r}\frac{\Phi_g(z)}s\,\dA(z)
\]
over polynomials $g$. For a fixed $r<1$, polynomial
approximation gives the same infimum over holomorphic $g$ on
$\D$. The following zero-packing rigidity estimate is immediate.

\begin{lem}
\label{lem:ajm-folding-rigidity}
For every $0<c_0<\rho_{\mathbb H}$ and fixed $b_0\in(0,1)$,
there exists $r_0<1$, independent of $f\in\mathcal S$, such that
\begin{equation}\label{eq:ajm-absolute-mass}
 P_r+N_r\ge c_0I_r-\mathrm O_{b_0}(1),
 \qquad r_0<r<1.
\end{equation}
The implicit constant can be chosen independently of $f$ and $r$.
\end{lem}
\begin{proof}
For $x\ge0$, $|x^2-1|\ge(x-1)^2$. Consequently,
\[
 |j_f(z)|\ge\Phi_{g_f}(z),
 \qquad g_f=z\Nop(f)\in\operatorname{Hol}(\D).
\]
Since $|z|<1$ on $\mathcal A_r$,
\[
 P_r+N_r\ge\int_{\mathcal A_r}\frac{\Phi_{g_f}(z)}s\,\dA(z).
\]
The part with $|z|\le b_0$ is bounded uniformly over
$f\in\mathcal S$, by the standard pre-Schwarzian distortion
bound $s|\Nop(f)|\le6$. The defining lower limit for
$\rho_{\mathbb H}$ now gives, uniformly for $r$ sufficiently
close to $1$,
\[
 P_r+N_r\ge
 (\rho_{\mathbb H}-o(1))L(r)-\mathrm O_{b_0}(1).
\]
As $I_r=L(r)+\mathrm O_{b_0}(1)$, this proves the assertion.
\end{proof}

The lemma controls absolute Jacobian mass, not its sign. We now establish a uniform height gap at saddles with
$j_f\ge\tau$, pair these saddles with maxima along complete target
fibers, and bound the complementary positive Jacobian mass.

We first establish a uniform local height gain at marked saddles.

\begin{lem}
\label{lem:buffered-saddle-height}
Fix $\tau>0$ and $b_0\in(0,1)$. There are positive constants
$h_\tau$, $\varepsilon_\tau$, and $K_\tau>1$, independent of
$f\in\mathcal S$ and $r$ with $r$ sufficiently close to $1$, with
the following property. If
\[
 b_0<|z_0|<r_*,\qquad
 1-r_*^2=K_\tau(1-r^2),\qquad
 j_f(z_0)=(1-|z_0|^2)^2|z_0\Nop(f)(z_0)|^2-1\ge\tau,
\]
then $z_0$ is a nondegenerate saddle of $u_{\xi,r}$ with
$\xi=G_r(z_0)$. The two opposite local upper sectors of this saddle each
contain a point of height at least $u_{\xi,r}(z_0)+h_\tau$
within Euclidean distance $\varepsilon_\tau(1-|z_0|^2)$ of $z_0$.
The abstract saddle--maximum matching is given in
Lemma~\ref{lem:abstract-one-pole-pairing}; its application to complete
fibers is in Lemma~\ref{lem:buffered-weighted-area-pairing}.
\end{lem}
\begin{proof}
Put $s_0=1-|z_0|^2$ and write $z=z_0+s_0w$.
A fixed small disk $|w|\le\varepsilon_0$ remains in $\D$
for every $z_0\in\D$. After shrinking $\varepsilon_0$,
Koebe distortion bounds the normalized holomorphic function
\[
 F_{z_0}(w)=\frac{f(z_0+s_0w)-f(z_0)}{s_0 f'(z_0)}
\]
and its first three derivatives uniformly on this disk.
On the untapered region $G_r=f+sf'/\bar z$, so
\[
 f(z_0+s_0w)-\xi
 =s_0f'(z_0)\left(F_{z_0}(w)-\frac1{\bar z_0}\right).
\]
As $|z_0|\ge b_0$ and $|z_0|<1$, the expression inside the
parentheses is bounded away from zero for sufficiently small
fixed $\varepsilon_0$. Similarly,
\[
 1-|z_0+s_0w|^2
 =s_0\bigl(1-2\re(\bar z_0w)-s_0|w|^2\bigr),
\]
whose second factor is bounded away from zero. It follows that
the third derivatives of the normalized potential
$u_{\xi,r}(z_0+s_0w)-u_{\xi,r}(z_0)$ are universally bounded.
Choosing $K_\tau>1$ sufficiently large keeps this disk in the
untapered region $1-|z|^2>1-r^2$ whenever $s_0>K_\tau(1-r^2)$.

At the saddle, the Hessian in the $w$ coordinate has trace $-4$
and determinant $-4j_f(z_0)$. Indeed, these follow from
$\Delta u_{\xi,r}=-4/s^2$ and the exact critical-point
Jacobian--Hessian identity. Its positive eigenvalue is therefore
\[
 \lambda_+=2\bigl(\sqrt{1+j_f(z_0)}-1\bigr)
 \ge 2\bigl(\sqrt{1+\tau}-1\bigr)>0.
\]
For a uniform quantitative choice, fix an upper bound $C_3$ for
all third derivatives of this normalized potential on
$|w|\le\varepsilon_0$. Put
\[
 \lambda_\tau=2(\sqrt{1+\tau}-1),\qquad
 t_\tau=\min\left\{\frac{\varepsilon_0}{2},
 \frac{\lambda_\tau}{2(1+C_3)}\right\},\qquad
 h_\tau=\frac{\lambda_\tau t_\tau^2}{4}.
\]
For a unit positive Hessian eigenvector $e$, Taylor's theorem gives,
for all $0<t\le t_\tau$,
\[
 u_{\xi,r}(z_0\pm s_0t e)-u_{\xi,r}(z_0)
 \ge\frac{\lambda_\tau t^2}{2}-\frac{C_3t^3}{6}
 \ge \frac{\lambda_\tau t^2}{4}>0.
\]
Thus both opposite rays stay strictly above the saddle level,
and their endpoints have height at least $h_\tau$ above it.
In particular, both endpoints lie above the intermediate level
$u_{\xi,r}(z_0)+h_\tau/2$.
Choose $K_\tau>(1-2\varepsilon_0-\varepsilon_0^2)^{-1}$.
Then the entire disk $|w|\le\varepsilon_0$, including both rays,
stays in the untapered region whenever $s_0>K_\tau(1-r^2)$.
The rays occupy the two opposite local upper sectors of the
nondegenerate saddle. This proves the claimed local estimates.
The complete-fiber pairing with distinct maxima follows from
Lemma~\ref{lem:abstract-one-pole-pairing}.
\end{proof}

For the global pairing, we may discard a planar-null set of targets
with coinciding critical heights.

\begin{lem}
\label{lem:generic-critical-heights}
Suppose $f$ is univalent across $\overline\D$ and $\Psi$ satisfies
 the hypotheses of the one-pole folding lemma. Outside a planar-null
 set of targets $\xi$, the critical points of $u_\xi$ are
 nondegenerate and have pairwise distinct finite critical heights.
\end{lem}
\begin{proof}
The critical points are precisely the inverse images of $\xi$
 under the smooth critical-target map $G$, away from the origin.
 Sard's theorem excludes a null set of critical values of $G$.
 For every other target outside $f(\partial\D)$, the fiber is finite
 and its inverse branches $z_j(\xi)$ are smooth locally.
 Write $H_j(\xi)=u_\xi(z_j(\xi))$. Since $z_j(\xi)$ is a critical
 point of $u_\xi$, the chain rule gives the envelope derivative
 \[
 \partial_\xi H_j(\xi)=\frac1{f(z_j(\xi))-\xi}.
 \]
 For $j\ne\ell$ these derivatives are different because $f$ is
 injective and $z_j\ne z_\ell$. Therefore the real-valued function
 $H_j-H_\ell$ has nowhere-vanishing differential on its domain.
 Its zero set is locally a smooth real curve and hence has planar
 measure zero. A countable collection of inverse-branch charts and
 the finitely many pairs on each chart suffice to cover all possible
 coincidences. Discarding their union proves the claim.
\end{proof}

\Needspace{5\baselineskip}
We first isolate the topological mechanism behind the quantitative folding argument. It applies to general strictly superharmonic potentials with at most one logarithmic pole.

\begin{lem}
\label{lem:abstract-one-pole-pairing}
Let $u$ be real-valued and $C^3$ on $\D\setminus\{p\}$,
strictly superharmonic there, and tending to $-\infty$ at
$\partial\D$.  There may be no exceptional point $p$, or
there may be a single logarithmic pole
$u(z)=-2\log|z-p|+v(z)$, where $v$ extends as a $C^3$
function across $p$.
Suppose the ordinary critical points are finite in number,
nondegenerate, and have pairwise distinct heights.  Let
$\mathcal M$ and $\mathcal S$ denote the ordinary maxima and
saddles, respectively, and let $\mathcal S_*\subset\mathcal S$.
Suppose there is $h>0$ such that, from each of the two local
upper sectors at every $q\in\mathcal S_*$, a path runs from
$q$ through points of height strictly above $u(q)$ to a point
of height at least $u(q)+h$.
Then there is an injection $q\mapsto m(q)$ from $\mathcal S$
into $\mathcal M$ with $u(m(q))>u(q)$ for every saddle and
\begin{equation}\label{eq:abstract-pairing-height}
 u(m(q))\ge u(q)+h,\qquad q\in\mathcal S_*.
\end{equation}
Consequently,
\begin{equation}\label{eq:abstract-pairing-weighted}
 \sum_{m\in\mathcal M}\e^{u(m)}
 -\sum_{q\in\mathcal S}\e^{u(q)}
 \ge(\e^h-1)\sum_{q\in\mathcal S_*}\e^{u(q)}.
\end{equation}
The possible pole is not counted as an ordinary maximum.
\end{lem}
\begin{proof}
For any finite regular height $t$, the superlevel set
$\{u>t\}$ is relatively compact in $\D$ after including the
pole (if present) at height $+\infty$.  A bounded
complementary hole would force an interior minimum of the
strictly superharmonic function $u$.  Thus each superlevel
component is simply connected.  As the level decreases,
ordinary maxima create components.  At every ordinary saddle,
the two upper sectors must lie in distinct components just
above the saddle height: if they belonged to the same
component, the subsequent one-handle attachment would
create a hole.  Hence every saddle merges two components.

Give each component the label of its highest birth maximum,
or the pole label of height $+\infty$ if it contains the pole.
At a merger retain the higher label and assign the lower label
to that saddle.  Labels are discarded only once.  The pole
label can never be discarded, so this defines an injection
from all saddles to ordinary maxima of strictly greater height.

If $q\in\mathcal S_*$, choose a regular level just above
$u(q)$ and below $u(q)+h$, with no intervening critical
height.  By the path hypothesis, each of the two components
which merge at $q$ contains a point where $u\ge u(q)+h$.
Its birth label therefore has height at least $u(q)+h$.
This holds for the lower label, which is assigned to $q$,
and proves \eqref{eq:abstract-pairing-height}.
Pairing each saddle with its assigned maximum and leaving
any remaining maxima unpaired gives
\eqref{eq:abstract-pairing-weighted}.
\end{proof}

For the complete potentials $u_{\xi,r}$, the abstract pairing lemma
applies to almost every target. Integrating its quantitative surplus
gives the required estimate.

\begin{lem}
\label{lem:buffered-weighted-area-pairing}
Let $\eta$ be the nonnegative compactly supported target cutoff
of Subsection~2.2, and define
\[
 \mathcal D_r=-\int_{\D}\eta(G_r(z))
  \e^{u_{G_r(z),r}(z)}J_{G_r}(z)\,\dA(z).
\]
For the buffered annulus
$A_*:=\{b_0<|z|<r_*\}$ of
Lemma~\ref{lem:buffered-saddle-height}, let
\[
 P_*^{(\tau)}=\int_{A_*\cap\{j_f\ge\tau\}}
                  \frac{j_f}{s|z|^2}\,\dA.
\]
Then, uniformly over $f\in\mathcal S_M$ extending univalently
past the closed disk,
\begin{equation}\label{eq:weighted-fiber-pairing}
 \mathcal D_r\ge(\e^{h_\tau/2}-1)P_*^{(\tau)}.
\end{equation}
\end{lem}
\begin{proof}
Exclude the planar-null set consisting of $f(\partial\D)$, the
critical values of $G_r$, and targets with coincident critical heights
(Lemma~\ref{lem:generic-critical-heights}). For every remaining target
$\xi$, the fiber $G_r^{-1}(\xi)$ is finite: $G_r(z)\to\infty$ as
$z\to0$, while $G_r(z)-f(z)\to0$ as $|z|\to1$.

The potential $u_{\xi,r}$ satisfies
Lemma~\ref{lem:abstract-one-pole-pairing}; at every marked saddle,
Lemma~\ref{lem:buffered-saddle-height} supplies its two upper paths.
Applying the abstract lemma with the reduced gap $h_\tau/2$ gives
\[
 \sum_{\mathrm{max}}\e^{u_{\xi,r}}-
 \sum_{\mathrm{saddle}}\e^{u_{\xi,r}}
 \ge (\e^{h_\tau/2}-1)
 \sum_{\substack{\mathrm{marked}\\\mathrm{saddles}}}\e^{u_{\xi,r}}.
\]
The possible logarithmic pole is not counted as an ordinary maximum.

For the integrals, apply the nonnegative area formula to each Jacobian
sign on compact source exhaustions. The target cutoff and the quadratic
outer taper give absolute integrability by Subsection~2.2, so the
exhaustion limits may be subtracted and the null exceptional set ignored.
Integrating the preceding fiber inequality against $\eta(\xi)$ yields
\[
 \mathcal D_r\ge(\e^{h_\tau/2}-1)
 \int_{A_*\cap\{j_f\ge\tau\}}\eta(G_r)
 \e^{u_{G_r,r}}J_{G_r}\,\dA.
\]
On $A_*$, $\eta(G_r)=1$ and
$\e^{u_{G_r,r}}J_{G_r}=j_f/(s|z|^2)$, proving
\eqref{eq:weighted-fiber-pairing}.
\end{proof}

We also need a uniform bound for the complementary source region.

\begin{lem}
\label{lem:buffered-complete-fiber}
With the preceding notation put
\[
 I_*:=\int_{A_*}\frac{\dA}{s|z|^2},\qquad
 P_*:=\int_{A_*}\frac{(j_f)_+}{s|z|^2}\,\dA,
 \qquad N_*:=\int_{A_*}\frac{(j_f)_-}{s|z|^2}\,\dA.
\]
There is a constant $C_{b_0,\tau}$, independent of $r$, $M$,
and $f\in\mathcal S_M$, such that
\begin{equation}\label{eq:buffered-complete-estimates}
 \mathcal D_r\ge N_*-P_*-C_{b_0,\tau}.
\end{equation}
\end{lem}
\begin{proof}
We estimate the positive part of the weighted Jacobian in the
complement of $A_*$. On $A_*$ the target cutoff equals one and
$\mathcal D_r$ contributes exactly $N_*-P_*$. Since the negative
part of the Jacobian in the complementary region helps the
inequality, it suffices to bound
\[
 E:=\int_{\D\setminus A_*}\eta(G_r)
             \e^{u_{G_r,r}}(J_{G_r})_+\,\dA.
\]
In the central region $|z|\le b_0$, the taper is inactive.
The positive part is better behaved than the absolute Jacobian.
The classical pre-Schwarzian estimate $s|\Nop(f)|\le6$ gives
\[
 (j_f)_+=\bigl(s^2|z\Nop(f)|^2-1\bigr)_+
 \le36|z|^2.
\]
In particular $j_f\le0$ for $|z|\le 1/6$. Since $0\le\eta\le1$,
\[
 E_{\mathrm{center}}
 \le36\int_{|z|<b_0}\frac{\dA}{1-|z|^2}
 =36\log\frac1{1-b_0^2}=C_{b_0}.
\]
Thus the central positive mass is uniformly bounded without
using the target cutoff. The cutoff remains necessary for absolute
integrability of the negative mass at the central pole, and it
causes the separate $2\log M$ loss in the radial comparison.
For the remaining outer region, recall $s(r)=1-r^2$ and
write $L=\log(1/s)$. In the untapered buffer
$s(r)\le s\le K_\tau s(r)$, the same estimate
$|j_f|\le37$ gives
$E_{\mathrm{buffer}}\le C_{b_0}(1+\log K_\tau)$.
In the transition $s(r)/\e\le s\le s(r)$, the quantities
$B_r(L)$, $\frac{\diff B_r}{\diff L}(L)$, $A_r(s)/s$,
and $A_r'(s)$ are uniformly bounded,
while $A_r(s)\asymp s$. The exact folding-Jacobian formula and
the pre-Schwarzian distortion bound give
\[
 |\e^{u_{G_r,r}}J_{G_r}|
 \le\frac{C}{s}
 \qquad(s(r)/\e\le s\le s(r),\ |z|\ge b_0).
\]
The transition consequently contributes $O(1)$.
Finally, if $0<s\le s(r)/\e$, then $B_r(L)=2$,
$\frac{\diff B_r}{\diff L}(L)=0$, $A_r(s)=s/2$,
and $W_r(s)=\exp(1+Q(1))s^2/s(r)$.
The exact Jacobian remainder is now
$R_r(s)=1/(1-s)+2/s$, and the mixed coefficient equals $2/s$.
It follows once more from $s|\Nop(f)|\le6$ that
\[
 |\e^{u_{G_r,r}}J_{G_r}|
 \le\frac{C}{s(r)}
 \qquad(0<s\le s(r)/\e,\ |z|\ge b_0).
\]
The hyperbolic outer tail has Euclidean radial thickness
$O(s(r))$, so its contribution is also $O(1)$.
Summing these four estimates gives
$E\le C_{b_0,\tau}$, as desired.
\end{proof}

\subsection{The Littlewood--Paley comparison}
The endpoint $j=0$ of the differentiation rule discussed in
\cite{DebnathHedenmalmDirichlet} follows from the full-disk
Littlewood--Paley identity and a uniform Bloch remainder estimate.
For $f\in\mathcal S$, put $g=\log f'$ and use the ordinary
dilation $g_r(z)=g(rz)$.  Thus $g_r'(z)=r\Nop(f)(rz)$.
Since $g_r(0)=0$, Green's formula on the unit disk gives
\[
 \|g_r\|_{H^2}^2
 =\int_\D |g_r'(z)|^2\log\frac1{|z|^2}\,\dA(z).
\]
The universal pre-Schwarzian distortion bound yields
\[
 (1-|z|^2)|g_r'(z)|\le6.
\]
Writing $s=|z|^2$, the difference of the logarithmic and
linear kernels is $-\log s-(1-s)\ge0$.  Therefore
\[
 \begin{aligned}
 0&\le \|g_r\|_{H^2}^2
 -\int_\D(1-|z|^2)|g_r'(z)|^2\,\dA(z)\\
 &\le36\int_\D
 \frac{\log(1/|z|^2)-(1-|z|^2)}{(1-|z|^2)^2}
 \,\dA(z)=:C_{\mathrm{LP}}=36.
 \end{aligned}
\]
Indeed, using $u=|z|^2$, normalized area measure, and the power
series for $-\log(1-v)$, we find
\[
 \int_0^1\frac{-\log u-(1-u)}{(1-u)^2}\,\diff u
 =\sum_{n=2}^{\infty}\frac1{n(n-1)}=1.
\]
Changing variables $w=rz$ gives
\[
 \begin{aligned}
 \|(\log f')_r\|_{H^2}^2
 &\le\int_{|w|<r}
 \left(1-\frac{|w|^2}{r^2}\right)|\Nop(f)(w)|^2\,\dA(w)
 +C_{\mathrm{LP}}\\
 &\le\int_{|w|<r}(1-|w|^2)|\Nop(f)(w)|^2\,\dA(w)
 +C_{\mathrm{LP}}.
 \end{aligned}
\]
Thus the Littlewood--Paley comparison transfers any bounded-image
strict nonlinearity-energy estimate to the circular $H^2$ norm,
with an absolute additive error.

\subsection{The cylindrical hat trick and the unrestricted class}
The geometric idea behind the cylindrical root transform is to narrow
 the channels through which the conformal image reaches infinity.
Passing to an $N$-fold cylindrical covering reduces the logarithmic
extent of these channels by a factor $N$, and hence diminishes the
cost of the target cutoff in the folding estimate. Remarkably, this
modification leaves asymptotic variance unchanged, since the
corresponding logarithmic derivatives differ by a lower-order term.
The construction thus reduces the geometric cost of access to infinity
without sacrificing the boundary oscillation that we seek to estimate.

For $f\in\mathcal S$ and a positive integer $N$, define its
\emph{cylindrical root transform} by
\[
 \mathcal C_N f(z)=z\left(\frac{f(z^N)}{z^N}\right)^{1/N},
\]
where the holomorphic root equals $1$ at $0$. This is a normalized
schlicht map, and
\[
 (\mathcal C_Nf(z))^N=f(z^N),\qquad
 \|\mathcal C_Nf\|_\infty\le \|f\|_\infty^{1/N}
\]
whenever $f$ is bounded. Put
\[
 u_f=\log f',\qquad v_f=\log(f/z),\qquad c_N=1-N^{-1},
\]
with the logarithms normalized to vanish at $0$. Direct differentiation
shows that
\begin{equation}\label{eq:cylindrical-identity}
 \log(\mathcal C_Nf)'(z)
 =\bigl(u_f-c_Nv_f\bigr)(z^N).
\end{equation}
By the theorem of Duren and Leung \cite[Theorem~4]{DurenLeung1979},
the logarithmic coefficients of normalized schlicht functions satisfy
the sharp bound
\begin{equation}\label{eq:cylindrical-H2}
 \sup_{f\in\mathcal S}\|\log(f/z)\|_{H^2}
 \le C_0=\pi\sqrt{\frac23}.
\end{equation}
This bound is attained by the Koebe function.
In particular, for $0<t<1$, the triangle inequality gives
\begin{equation}\label{eq:cylindrical-norm}
 \left|\|(u_f)_{t^N}\|_{H^2}
 -\|(\log(\mathcal C_Nf)')_t\|_{H^2}\right|
 \le C_0.
\end{equation}

\subsection{Strict finite-scale deduction and proof of the main theorem}
\label{subsec:strict-deduction}

\begin{prop}
\label{prop:strict-finite-scale-deduction}
Fix $0<c_0<\rho_{\mathbb H}$ and $0<\tau<c_0/2$, and let
$h_\tau>0$ be as in Lemma~\ref{lem:buffered-saddle-height}. Set
\[
 \delta_\tau=(c_0-2\tau)\tanh\frac{h_\tau}{4}>0.
\]
For every $M\ge1$,
\begin{equation}\label{eq:strict-bounded-conditional}
 \int_{|z|<r}(1-|z|^2)|\Nop(f)(z)|^2\,\dA(z)
 \le (1-\delta_\tau)L(r)+C_\tau(1+\log M),
 \qquad f\in\mathcal S_M,
\end{equation}
where the lower admissible radius and $C_\tau$ are independent
of $M$ and $f$.
The corresponding unrestricted finite-radius bound is
\begin{equation}\label{eq:strict-unrestricted-conditional}
 \sup_{f\in\mathcal S}\|(\log f')_r\|_{H^2}^2
 \le(1-\delta_\tau)L(r)+O_\tau(\sqrt{L(r)}).
\end{equation}
In particular, these estimates imply
$\Sigma_{\max}^2\le1-\delta_\tau<1$.
\end{prop}
\begin{proof}
Work first with $f\in\mathcal S_M$ extending univalently past the
closed disk. Replace $r$ by the buffered radius $r_*$ in
Lemma~\ref{lem:ajm-folding-rigidity}. Since
$L(r_*)=L(r)-\log K_\tau$, its conclusion gives
\[
 P_*+N_*\ge c_0I_*-C_\tau,
 \qquad I_*=L(r)-\log K_\tau+O_{b_0}(1).
\]
Put $a_\tau=\exp(h_\tau/2)-1$. The complementary source estimate in
Lemma~\ref{lem:buffered-complete-fiber}
then supplies
\[
 \mathcal D_r\ge a_\tau P_*^{(\tau)},\qquad
 \mathcal D_r\ge N_*-P_*-C_\tau.
\]
Since $P_*\le\tau I_*+P_*^{(\tau)}$ and
$N_*-P_*=(P_*+N_*)-2P_*$, we obtain
\[
 \left(1+\frac2{a_\tau}\right)\mathcal D_r
 \ge(c_0-2\tau)I_*-C_\tau.
\]
Equivalently,
\begin{equation}\label{eq:strict-defect-audited}
 \mathcal D_r\ge \delta_\tau L(r)-C_\tau,
 \qquad \frac{a_\tau}{a_\tau+2}=\tanh\frac{h_\tau}{4}.
\end{equation}
The signed Jacobian identity and the vanishing of its integrated
mixed term give
\[
 \int_\D\eta(G_r)W_r(s)|\Nop(f)|^2\,\dA
 +\mathcal D_r
 =\int_\D\eta(G_r)W_r(s)R_r(s)\,\dA.
\]
Subsection~2.2 bounds the right side by $L(r)+2\log M+\Ordo(1)$.
Since $W_r=s$ on $|z|<r$, the cutoff is one for $|z|\ge b_0$,
and the omitted central nonlinearity energy is $\Ordo_{b_0}(1)$,
we obtain \eqref{eq:strict-bounded-conditional} for smooth $f$.
For a general $f\in\mathcal S_M$, use
$\rho^{-1}f_\rho$, with $\rho\to1^-$.
The admissible radius and constants are independent of $\rho$,
and $\|\rho^{-1}f_\rho\|_\infty\le M/\rho$, so the same estimate
passes to the limit.

We make the transfer to the unrestricted schlicht class explicit.
The Littlewood--Paley comparison, now with $C_{\mathrm{LP}}=36$,
first gives, for all $g\in\mathcal S_M$ and all $q$ sufficiently
close to $1$,
\begin{equation}\label{eq:strict-H2-bounded-transfer}
 \|(\log g')_q\|_{H^2}^2
 \le (1-\delta_\tau)L(q)+C_\tau(1+\log M)+36.
\end{equation}
For a positive integer $N$, let $F=\mathcal C_Ng$ and set
$q=r^{1/N}$. The root transform belongs to
$\mathcal S_{M^{1/N}}$, and $q\ge r$. The identity
\eqref{eq:cylindrical-identity} and the Duren--Leung estimate
\eqref{eq:cylindrical-H2} imply
\[
 \|(\log g')_r\|_{H^2}
 \le\|(\log F')_q\|_{H^2}+C_0,
 \qquad C_0=\pi\sqrt{\frac23}.
\]
Consequently, by \eqref{eq:strict-H2-bounded-transfer},
\begin{equation}\label{eq:strict-root-transfer}
 \|(\log g')_r\|_{H^2}
 \le\left((1-\delta_\tau)L(q)
  +C_\tau\left(1+\frac{\log M}{N}\right)+36\right)^{1/2}
  +C_0.
\end{equation}
Since
\[
 1-r^{2/N}\ge\frac{1-r^2}{N},
 \qquad L(q)\le L(r)+\log N,
\]
we may choose $N=\max\{1,\lceil\log M\rceil\}$. Then
$\log M/N\le1$ and $\log N\le\log(1+\log M)$.
Squaring \eqref{eq:strict-root-transfer} gives, uniformly
in $M$ and $g\in\mathcal S_M$,
\begin{equation}\label{eq:strict-root-final}
 \|(\log g')_r\|_{H^2}^2
 \le (1-\delta_\tau)L(r)
 +O_\tau\bigl(\sqrt{L(r)}+\log(2+\log M)+1\bigr).
\end{equation}
The error term is lower order in $L(r)$, so the leading coefficient is unchanged.

Finally, let $f\in\mathcal S$ be arbitrary, put $t=\sqrt r$, and
consider $g=t^{-1}f_t$. Then $g\in\mathcal S_{M_t}$ with
\[
 M_t\le\frac1{(1-t)^2},\qquad
 \log g'(t\zeta)=\log f'(r\zeta).
\]
The Koebe growth bound and the identities
\[
 L(t)=L(r)+\log(1+r),\qquad
 \log M_t\le 2\log\frac1{1-t}
 \le 2L(r)+2\log4
\]
show that $\log(2+\log M_t)=O(\log(2+L(r)))$
as $r\to1^-$. Apply \eqref{eq:strict-root-final} to $g$ at radius
$t$. Since $\log(2+L(r))=O(\sqrt{L(r)})$, this proves
\eqref{eq:strict-unrestricted-conditional}, with constants and
an admissible radius independent of $f$. Dividing by $L(r)$ and
letting $r\to1^-$ yields the stated variance bound.
\end{proof}

\Needspace{5\baselineskip}
We now prove the main finite-scale theorem.

\begin{proof}
Fix constants $0<c_0<\rho_{\mathbb H}$ and $0<\tau<c_0/2$.
Proposition~\ref{prop:strict-finite-scale-deduction}, together with
its uniform Littlewood--Paley and cylindrical-root transfer, gives
\eqref{eq:strict-main-intro} with
$\Sigma_{\mathrm{fold}}^2=1-\delta_\tau$, where $\delta_\tau>0$.
The choices of $c_0$ and $\tau$ are absolute, so the error term
and the admissible radius are uniform over $f\in\mathcal S$.
\end{proof}

The asymptotic-variance corollary follows from
\eqref{eq:strict-unrestricted-conditional} upon division by $L(r)$
and letting $r\to1^-$. 

The infinitesimal zero-packing identity also restricts the size of any
possible strict improvement.

\begin{rem}
\label{rem:zero-packing-upper-constraint}
Denote by
$\Sigma_{\mathrm{inf}}^2$ the universal infinitesimal asymptotic
variance of the Beurling transform. The zero-packing identity of
\cite{HedenmalmAJM}, together with the ordinary/tight discrepancy
equality of~\cite{WennmanDiscrepancy}, asserts that
\[
 \Sigma_{\mathrm{inf}}^2=1-\rho_{\mathbb H}.
\]
The unrestricted maximal variance is at least this infinitesimal
constant:
\[
 \Sigma_{\max}^2\ge\Sigma_{\mathrm{inf}}^2.
\]
Indeed, for each $0<\kappa<1$, the quasiconformal subclass and
the monotonicity quoted above give
\[
 \Sigma_{\max}^2\ge \Sigma^2(\kappa)
 \ge \kappa^2\Sigma_{\mathrm{inf}}^2.
\]
The first comparison uses the standard equivalence of leading
logarithmic variance under the interior and exterior conformal
normalizations. The second follows from Ivrii's monotonicity
\cite[Theorem~1.4]{IvriiQuasicircles}.
Letting $\kappa\to1^-$ gives the asserted comparison.
Taking $\delta=1-\Sigma_{\mathrm{fold}}^2>0$ from
Theorem~\ref{thm:main-finite-scale}, we obtain
\begin{equation}\label{eq:zero-packing-upper-constraint}
 0<\delta\le1-\Sigma_{\max}^2\le\rho_{\mathbb H}\le\rho_{12,3}\simeq 0.0523107153.
\end{equation}
The final estimate is the fundamental-cell upper bound proved in
Proposition~\ref{prop:12-3-cell-upper-bound} below.
Thus zero packing imposes an \emph{upper} restriction on the
variance improvement supplied by folding, whereas
Lemma~\ref{lem:ajm-folding-rigidity} supplies a \emph{lower}
rigidity bound for the absolute Jacobian mass. In particular,
the coefficient $\delta_\tau$ of
Proposition~\ref{prop:strict-finite-scale-deduction} is compatible with
\eqref{eq:zero-packing-upper-constraint}, since
$\delta_\tau<c_0<\rho_{\mathbb H}$.
\end{rem}

\section{Hyperbolic zero packing, Neumann cells, and candidate domains}
\label{sec:extremal-domain-program}

We distinguish the periodic $\{12,3\}$ benchmark, a conjectural
Neumann-cell lower bound, and exploratory boundary domains. None is
used to prove the strict folding theorem of Section~2.

\subsection{Three variance quantities}
\label{subsec:three-variance-problems}

It is useful to distinguish the infinitesimal constant, the
variance at prescribed quasiconformal dilatation,
and the unrestricted maximal variance.
Let $\Sigma\langle\kappa\rangle$ denote normalized exterior conformal
maps $\psi(z)=z+\mathrm O(1)$ at infinity which admit an extension
with complex dilatation bounded in absolute value by $\kappa$.
For $0<\kappa<1$ set
\[
 \Sigma^2(\kappa):=\sup_{\psi\in\Sigma\langle\kappa\rangle}
 \limsup_{R\to1^+}
 \frac{\displaystyle\int_{\T}|\log\psi'(R\zeta)|^2\,\diff s(\zeta)}
 {\displaystyle\log\frac{R^2}{R^2-1}}.
\]
The corresponding infinitesimal constant is
\[
 \Sigma_{\mathrm{inf}}^2
 :=\lim_{\kappa\to0^+}\frac{\Sigma^2(\kappa)}{\kappa^2}
 =1-\rho_{\mathbb H}<1.
\]
The last identity follows from the geometric zero-packing
relation of~\cite{HedenmalmAJM} and the equality of the ordinary
and tight hyperbolic discrepancy densities in~\cite{WennmanDiscrepancy}.  The asymptotic behavior and monotonicity
properties of $\Sigma^2(\kappa)$ are discussed in
\cite{IvriiQuasicircles}.  The constant $\Sigma_{\max}^2$ of the
present paper instead allows arbitrary normalized schlicht maps.
The strict upper bound for this constant does not identify
$\lim_{\kappa\to1^-}\Sigma^2(\kappa)$.
The finite-scale theorem of this paper does not require a separate
periodic-spectrum argument.

\subsection{The regular hyperbolic $\{12,3\}$ zero lattice}
\label{subsec:12-3-zero-lattice}

We next describe an especially structured benchmark for the
\emph{infinitesimal} extremal problem.  Consider the regular
hyperbolic tessellation $\{12,3\}$ in the Poincar\'e disk.
Three dodecagons meet at every vertex.  We consistently use the
normalized hyperbolic area measure of~\cite{HedenmalmAJM},
\[
 \diff A_{\mathbb H}(z)
  :=\frac{\diff A(z)}{(1-|z|^2)^2},\qquad
 \diff A(z):=\frac{\diff x\,\diff y}{\pi}.
\]
The curvature $-1$ area element is
$\diff A_{-1}=4\pi\,\diff A_{\mathbb H}$.
Gauss--Bonnet therefore gives normalized area
\[
 \operatorname{Area}_{\mathbb H}(\mathcal F)
 =\frac1{4\pi}\left((12-2)\pi-12\frac{2\pi}{3}\right)
 =\frac12,
\]
so its centers have density $2$ per unit normalized hyperbolic
area.  The orientation-preserving
symmetry group is the triangle group
$\Gamma=\Delta(2,3,12)$.
The lattice is a single orbit $\Gamma z_0$ of an order-twelve
elliptic point.

Let $J$ be a Hauptmodul on the quotient orbifold, normalized to
take the values $0$, $1$, and $\infty$ at elliptic orbits of
orders $2$, $3$, and $12$, respectively.  With compatible roots on
the simply connected disk, put
\begin{equation}\label{eq:12-3-differential}
 f_{12}(z)=C\frac{J'(z)}{J(z)^{1/2}(J(z)-1)^{2/3}}.
\end{equation}
The local branching orders show that $f_{12}$ has simple zeros
at the order-twelve elliptic points and no other zeros.
The differential $f_{12}(z)\,\diff z$ transforms with a unitary
character of $\Gamma$.  On an appropriate finite cover it is a
holomorphic section of the canonical bundle.  A convenient
algebraic model for a four-cell genus-three quotient is
\[
 X:\quad u^4=t^3+1,\qquad
 \omega=\frac{\diff t}{u^2}.
\]
It realizes a divisor of four simple zeros on the compact surface.
The quotient construction is valuable because the analytic
constraint on the zeros is imposed globally, rather than cell by cell.

\begin{figure}[t]
 \centering
 \includegraphics[width=.72\linewidth]{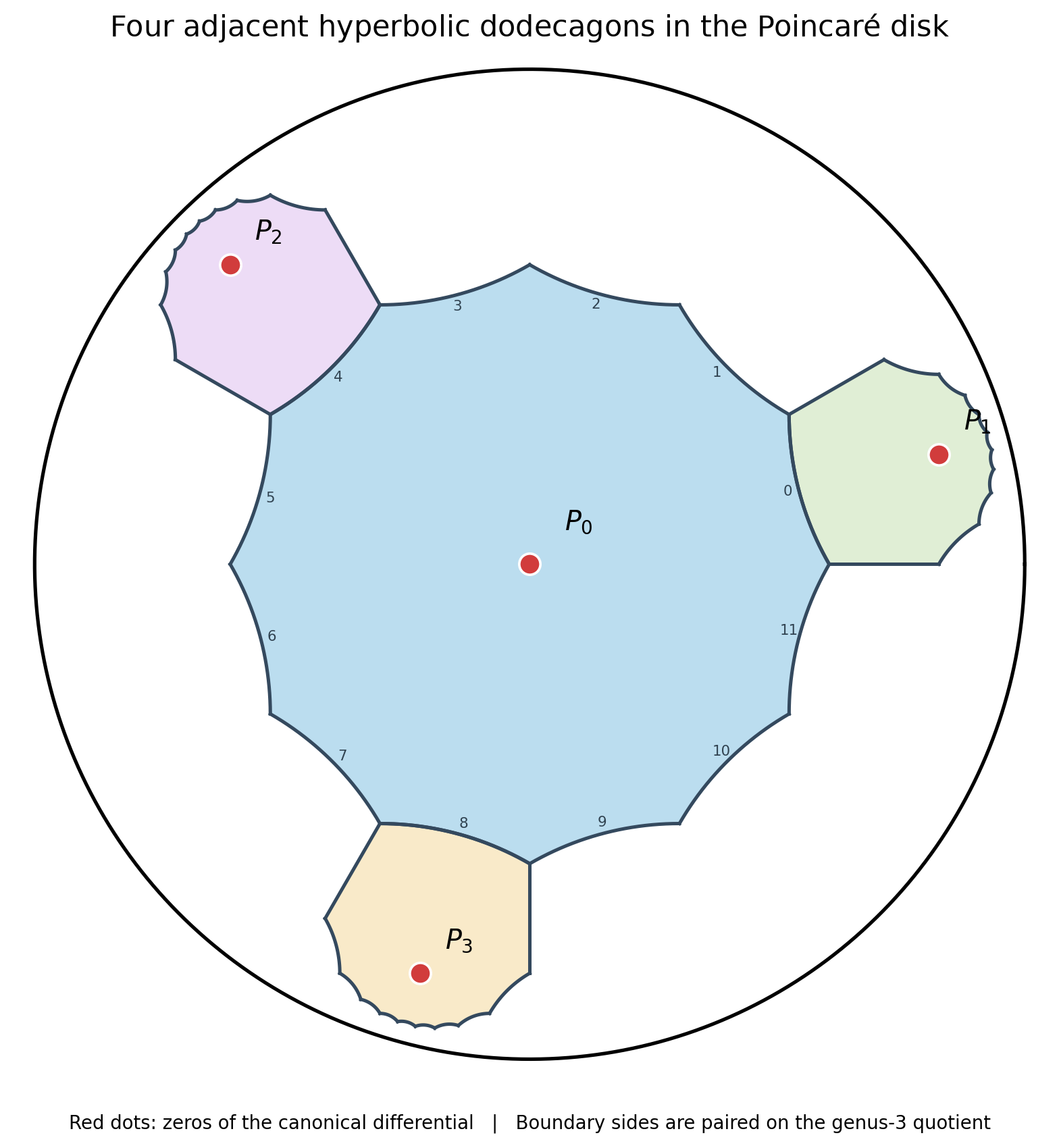}
 \caption{Four adjacent cells of the hyperbolic $\{12,3\}$ tessellation
 with marked centers.  Suitable boundary identifications produce a
 genus-three quotient.  The picture specifies a periodic zero-lattice
 candidate, not an established global optimizer.}
 \label{fig:12-3-four-cells}
\end{figure}

Let $\mathcal F$ be a regular hyperbolic dodecagon of normalized area $\frac12$,
and write $\langle\cdot\rangle_{\mathcal F}$ for its normalized
hyperbolic-area average.  Put
\[
 H_0(z)=(1-|z|^2)|f_{12}(z)|,
 \qquad \rho_{12,3}:=1-
 \frac{\langle H_0\rangle_{\mathcal F}^{\,2}}
      {\langle H_0^2\rangle_{\mathcal F}}.
\]
This quotient is independent of the multiplicative normalization of
$f_{12}$.

The regular cell gives an explicit upper bound for the hyperbolic
zero-packing discrepancy.

\begin{prop}
\label{prop:12-3-cell-upper-bound}
The hyperbolic zero-packing discrepancy satisfies
\begin{equation}\label{eq:12-3-cell-upper}
 0<\rho_{\mathbb H}\le\rho_{12,3}
 =1-\frac{\left(\int_{\mathcal F}H_0\,\diff A_{\mathbb H}\right)^2}
 {\frac12\int_{\mathcal F}H_0^2\,\diff A_{\mathbb H}}.
\end{equation}
In particular, the cell integral is an admissible, explicitly
specified upper bound, irrespective of whether the cell arrangement
is globally optimal.
\end{prop}
\begin{proof}
Write $H_c=cH_0$ for $c>0$.  Its cell-averaged squared
discrepancy is
\[
 \langle(H_c-1)^2\rangle_{\mathcal F}
 =c^2\langle H_0^2\rangle_{\mathcal F}
  -2c\langle H_0\rangle_{\mathcal F}+1.
\]
The minimum is attained at
$c=\langle H_0\rangle_{\mathcal F}/
\langle H_0^2\rangle_{\mathcal F}$ and is the right side
of \eqref{eq:12-3-cell-upper}.  Since $H_0$ is
$\Gamma$-invariant, it descends to a continuous function on the
compact hyperbolic orbifold $\mathbb D/\Gamma$.  The hyperbolic
circle equidistribution theorem, followed by radial Ces\`aro
averaging in the logarithmic variable
$\log(1/(1-|z|^2))$, identifies its normalized large-disk
logarithmic averages with normalized orbifold-area averages.
The same applies to $H_0^2$.  One way to see the circle
statement is to use the spectral decomposition of the Laplacian
on the compact quotient: the spherical mean of a nonconstant
Laplacian eigenfunction decays to zero with the hyperbolic radius,
whereas the constant eigenfunction is unchanged.  Approximation
by smooth functions then gives convergence for the continuous,
$\Gamma$-invariant functions $H_0$ and $H_0^2$.  The logarithmic
radial averages inherit the same limit by Ces\`aro averaging.
Thus the disk competitors $cf_{12}$ give the stated upper bound
for $\rho_{\mathbb H}$.
\end{proof}

Numerical evaluation of the fundamental-cell integrals gives the
computed values
\begin{equation}\label{eq:12-3-numerical-discrepancy}
 \rho_{12,3}\simeq0.0523107153,
 \qquad 1-\rho_{12,3}\simeq0.9476892847.
\end{equation}
An a posteriori enclosure, based on the Neumann boundary residual
and interval quadrature, is given in the accompanying numerical supplement.
In particular, it establishes the rigorous inequality
$\rho_{12,3}<0.052314$. The displayed ten-digit decimal values
are numerical approximations, not ten-digit error enclosures.
For comparison, Conjecture~\ref{conj:neumann-cell-sharp}
proposes the lower benchmark
\[
 \rho_*=0.0522117855\ldots,\qquad
 \rho_{12,3}-\rho_*\simeq9.89298\times10^{-5}.
\]
This proximity motivates studying the structured lattice, but it
neither proves the conjectural lower bound nor makes the periodic
configuration an extremizer.
\FloatBarrier
The cell upper bound does not establish optimality of the periodic
configuration. In fact, a disk-variation argument communicated by
Oleg Ivrii suggests that the exact ratio of a fixed automorphic
one-form can be strictly improved when nonperiodic perturbations
are allowed. In particular, the regular lattice should be regarded
as a benchmark rather than a conjectured global minimizer.
A natural starting point for investigating possible deformations
of its zeros is the finite-disk
Euler--Lagrange equation of~\cite[(1.5.3)]{HedenmalmAJM}.
Let $\mathbf P_r$ denote the ordinary Bergman projection on
$\D_r=\{z:|z|<r\}$ and put
$\mu_F=F/|F|$ almost everywhere, with the value at zeros set to
zero.  For a finite-disk discrepancy minimizer $F^{(r)}$, the
equation reads
\[
 \mathbf P_r\bigl((1-|z|^2)F^{(r)}-\mu_{F^{(r)}}\bigr)=0.
\]
The projection of the first term is explicit, so equivalently
\begin{equation}\label{eq:finite-disk-EL}
 \mathbf P_r\mu_{F^{(r)}}(z)
 =(1-r^2)F^{(r)}(z)
 +\frac{r^2}{z^2}\int_0^z\zeta F^{(r)}(\zeta)\,\diff\zeta.
\end{equation}
If $r_j\to1^-$ and $F^{(r_j)}$ converges locally uniformly to a
nonzero holomorphic $F$, then the phase functions converge almost
everywhere, and bounded convergence in the Bergman kernels gives
the limiting necessary stationarity condition
\begin{equation}\label{eq:limiting-zero-packing-EL}
 \mathbf P\mu_F(z)=\frac1{z^2}\int_0^z
 \zeta F(\zeta)\,\diff\zeta,
 \qquad z\in\D,
\end{equation}
where $\mathbf P$ is the Bergman projection on $\D$.
Equivalently, by differentiation,
\begin{equation}\label{eq:limiting-zero-packing-fixed-point}
 F(z)=2\int_\D\frac{\mu_F(\zeta)}
 {(1-z\overline\zeta)^3}\,\diff A(\zeta).
\end{equation}
In particular, every zero $a$ of a nonzero limiting function
satisfies the vanishing of the integral on the right at $z=a$.
Neither local uniform compactness of the optimizers after suitable
normalization nor optimality of a solution of this equation is
asserted.  The equation motivates studying deformations of the
$\{12,3\}$ zero set and periodic approximations on larger covers.

\subsection{A Neumann-cell conjecture for the discrepancy}
\label{subsec:neumann-cell-conjecture}

A possible geometric route to a sharp lower bound for
$\rho_{\mathbb H}$ is suggested by gradient flow.  We continue to use
the \emph{normalized} hyperbolic area measure
$\diff A_{\mathbb H}=\diff A/(1-|z|^2)^2$.
For lengths and gradients we use the ordinary curvature $-1$
metric, with
$\diff s_{-1}=2|\diff z|/(1-|z|^2)$ and corresponding
gradient $\nabla_{-1}$.  Its area measure is
$\diff A_{-1}=4\pi\diff A_{\mathbb H}$, and its Laplace
operator is $\Delta_{\mathbb H}=(1-|z|^2)^2\partial\bar\partial$.
For a holomorphic one-form $f(z)\,\diff z$, put
\[
 u(z)=\log\bigl((1-|z|^2)|f(z)|\bigr).
\]
In the distributional sense,
\begin{equation}\label{eq:neumann-zero-charge}
 \Delta_{\mathbb H}u\,\diff A_{\mathbb H}
 =\frac12\sum_a m_a\delta_a-\diff A_{\mathbb H},
\end{equation}
where $m_a$ is the multiplicity of a zero. For simple zeros, a
regular basin of attraction of the downward gradient flow has
boundary arcs tangent to $\nabla_{-1}u$. Hence its normal
derivative vanishes, and the divergence theorem gives normalized
hyperbolic area $\frac12$ for a complete basin containing one simple
zero (and $m/2$ for one zero of multiplicity $m$).
This is the hyperbolic counterpart of the equal-area gradient-flow
transportation to random zeros studied by Nazarov, Sodin, and
Volberg~\cite{NazarovSodinVolbergGradient}. Such a basin decomposition
requires appropriate regularity and control of exceptional
trajectories. Its extension to arbitrary near-minimizers on the
hyperbolic disk, including multiple zeros and boundary effects,
remains to be justified.

\begin{figure}[tbp]
 \centering
 \includegraphics[width=.97\linewidth]{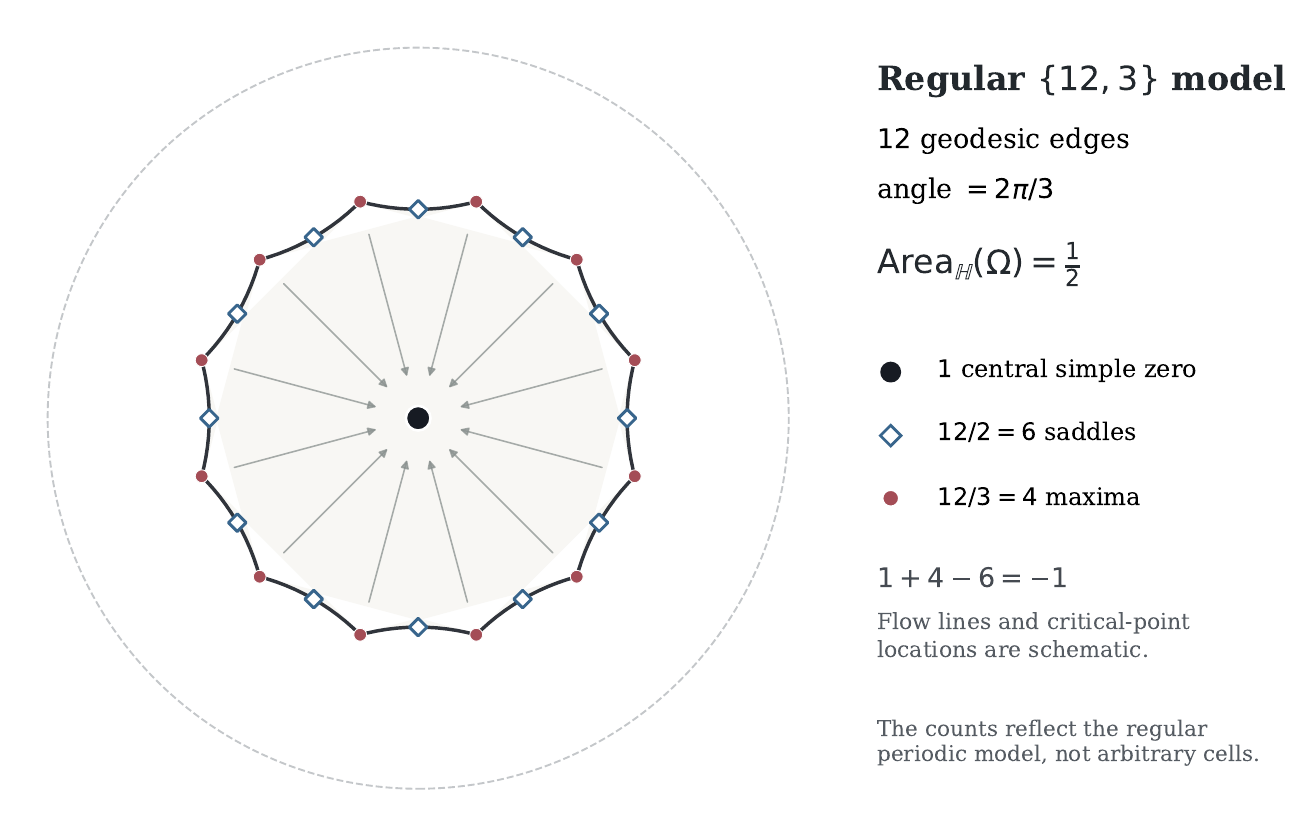}
 \caption{Schematic complete Neumann cell of normalized hyperbolic area $\frac12$
 in the regular $\{12,3\}$ model, shown inside the Poincar\'e disk.
 The central point is a simple zero, the open diamonds mark
 expected saddles on the twelve edges, and filled circles indicate
 expected local maxima at the twelve vertices.  Downward gradient
 curves head toward the zero; the cell sides are Neumann streamlines.
 Each edge saddle is shared by two cells and each vertex maximum
 by three. The effective counts per cell are therefore one zero,
 six saddles, and four maxima.  The flow lines and placement of
 critical points are illustrative, not a numerically solved gradient
 field.  A general Neumann cell may be irregular and need not have
 these critical-point counts.}
 \label{fig:neumann-cell-12-3}
\end{figure}

\paragraph{The geometry of one complete cell.}
Figure~\ref{fig:neumann-cell-12-3} depicts a single idealized Neumann
cell in the regular hyperbolic $\{12,3\}$ model.  The regular
dodecagon has twelve geodesic sides and interior angle
$2\pi/3$, so the Gauss--Bonnet formula, divided by
$4\pi$ to normalize the area, gives
\[
 \operatorname{Area}_{\mathbb H}(\mathcal F)
  =\frac{(12-2)\pi-12\cdot\frac{2\pi}{3}}{4\pi}
  =\frac12.
\]
The same area follows analytically for any complete Neumann basin
containing one simple zero: integrate
\eqref{eq:neumann-zero-charge} over the basin and use
$\partial_{\mathbf n}u=0$ on its boundary.  Thus the
normalized area $\frac12$ is not a consequence of dodecagonal symmetry, whereas
the depicted critical-point pattern is.

\FloatBarrier

The sharing factors match the index count on a compact quotient:
\[
 1+\frac{12}{3}-\frac{12}{2}=-1
 =-2\operatorname{Area}_{\mathbb H}(\mathcal F).
\]
On a compact hyperbolic surface this expresses the Morse relation
$N_{\mathrm{zero}}+N_{\mathrm{max}}-N_{\mathrm{saddle}}=\chi(X)$
for a generic gradient field, with each simple zero regarded as a
singular minimum.  The circular Neumann model discussed below is
rather different: its boundary is a degenerate critical level, so
its critical-point structure cannot be read as twelve saddles and
twelve maxima.  The conjectural inequality compares the
\emph{discrepancy} of cells, not their number of boundary critical
points.

For a sufficiently regular simply connected hyperbolic domain
$\Omega$ of normalized area $\frac12$ and a point $a\in\Omega$, let
$u_{\Omega,a}$ denote the solution, modulo additive constants, of
\begin{equation}\label{eq:neumann-cell-problem}
 \begin{cases}
 \Delta_{\mathbb H}u_{\Omega,a}\,\diff A_{\mathbb H}
    =\frac12\delta_a-\diff A_{\mathbb H}&\text{in }\Omega,\\
 \partial_{\mathbf n}u_{\Omega,a}=0&\text{on }\partial\Omega.
 \end{cases}
\end{equation}
Its scale-optimized discrepancy is
\begin{equation}\label{eq:neumann-cell-ratio}
 \mathcal R_N(\Omega,a)
 =1-\frac{\bigl(\int_\Omega \exp(u_{\Omega,a})\,
       \diff A_{\mathbb H}\bigr)^2}
 {\frac12\int_\Omega \exp(2u_{\Omega,a})\,
       \diff A_{\mathbb H}}.
\end{equation}
The expression is independent of the additive normalization of
$u_{\Omega,a}$ and is invariant under hyperbolic isometries.

The following conjecture proposes a sharp lower bound for the
scale-optimized discrepancy of an admissible hyperbolic Neumann cell.

\begin{conj}
\label{conj:neumann-cell-sharp}
For every admissible pair $(\Omega,a)$ as above,
\begin{equation}\label{eq:neumann-cell-lower}
 \mathcal R_N(\Omega,a)\ge\rho_*,\qquad
 \rho_*:=1-144\left(\operatorname{arctanh}\frac1{\sqrt3}
              -\frac1{\sqrt3}\right)^2
       =0.0522117855\ldots .
\end{equation}
Equality holds only for a hyperbolic disk of normalized area $\frac12$ centered
at $a$, up to a hyperbolic isometry.
\end{conj}

For the circular model, a M\"obius automorphism places the pole at
zero. Indeed, the normalized hyperbolic area of the Euclidean disk
of radius $R$ is $R^2/(1-R^2)$, which equals $1/2$ when
$R=3^{-1/2}$.  In these coordinates,
\[
 u_*(z)=\log\bigl(|z|(1-|z|^2)\bigr)+C,
 \qquad |z|<3^{-1/2}.
\]
Its normal derivative vanishes at the boundary because
$\partial_r\{r(1-r^2)\}=0$ at $r=3^{-1/2}$.
Optimizing the multiplicative constant yields the value $\rho_*$
in~\eqref{eq:neumann-cell-lower}. The circular value is an exact
model calculation, while the global comparison asserted in the
conjecture is open.

\Needspace{5\baselineskip}
The next proposition establishes positivity of the second
variation at the centered circular Neumann cell for smooth
area-preserving deformations with the zero fixed.

\begin{prop}
\label{prop:neumann-disk-local-minimum}
Let $R=1/\sqrt3$, and let $\Omega_\varepsilon$ be a smooth
one-parameter family of simply connected hyperbolic domains of
normalized area $\frac12$, with $\Omega_0=\{|z|<R\}$ and boundary
\[
 r=R+\varepsilon h(\theta)+\varepsilon^2k(\theta)
       +\ordo(\varepsilon^2),\qquad
 \int_0^{2\pi}h(\theta)\,\diff\theta=0.
\]
Keep the logarithmic pole fixed at $a=0$, and let
$u_\varepsilon$ solve the Neumann problem
\eqref{eq:neumann-cell-problem} on $\Omega_\varepsilon$.
For every nonzero real-valued smooth $h$ with zero mean,
\begin{equation}\label{eq:neumann-disk-second-variation}
 \mathcal R_N(\Omega_\varepsilon,0)
  =\rho_*+c[h]\varepsilon^2+\ordo(\varepsilon^2),
 \qquad c[h]>0.
\end{equation}
Thus the second variation is strictly positive along every
smooth area-preserving deformation with a nonzero first-order
boundary variation and the pole fixed at the center.  The
proposition does not assert a uniform local minimum in any
infinite-dimensional topology of domains.
\end{prop}
\begin{proof}
The circular potential is $u_0(z)=\log(|z|(1-|z|^2))$,
up to a constant.  Write $X=\e^{u_0}$ and
$H=\left.\partial_\varepsilon u_\varepsilon\right|_0$,
normalized by $H(0)=0$.  Since the pole and the curvature
source are fixed, $H$ is harmonic in $|z|<R$.
Differentiating the Neumann boundary condition gives
\[
 \partial_rH(R,\theta)=-u_0''(R)h(\theta)=9h(\theta).
\]
Writing $h(\theta)=\sum_{m\ne0}h_m\e^{\imag m\theta}$,
we therefore have
\begin{equation}\label{eq:neumann-first-shape-harmonic}
 H(r,\theta)=9R\sum_{m\ne0}\frac{h_m}{|m|}
      \left(\frac rR\right)^{|m|}\e^{\imag m\theta}.
\end{equation}
Put
\[
 M_p=\int_{|z|<R}X^p\,\diff A_{\mathbb H},\quad
 X_R=X(R),\quad \omega_R=\frac{R}{\pi(1-R^2)^2}
 \qquad(p=1,2).
\]
Let $M_p(\varepsilon)=\int_{\Omega_\varepsilon}
\e^{p u_\varepsilon}\,\diff A_{\mathbb H}$.  Write
$u_\varepsilon=u_0+\varepsilon H+\varepsilon^2K+
\ordo(\varepsilon^2)$ on compact subsets of the circular
cell, where $K$ is harmonic.  Differentiating the domain
integrals to second order, using $u_0'(R)=0$ and the
fixed-area constraint, yields
\[
\begin{aligned}
 [\varepsilon^2]M_p(\varepsilon)
 ={}&pK(0)M_p
 +\frac{p^2}{2}\int_{|z|<R}X^pH^2\,\diff A_{\mathbb H}\\
 &+p\omega_R X_R^p\int_0^{2\pi}h(\theta)H(R,\theta)
                  \,\diff\theta.
\end{aligned}
\]
Here the mean-value property evaluates the weighted integral
of the harmonic term $K$ as $K(0)M_p$.  Taking the
logarithmic ratio gives
\begin{equation}\label{eq:neumann-ratio-hessian}
\begin{aligned}
 [\varepsilon^2]\log\frac{M_2(\varepsilon)}{M_1(\varepsilon)^2}
 ={}&2\frac{\int_{|z|<R}X^2H^2\,\diff A_{\mathbb H}}{M_2}
  -\frac{\int_{|z|<R}XH^2\,\diff A_{\mathbb H}}{M_1}\\
 &+2\omega_R\left(\frac{X_R^2}{M_2}
                   -\frac{X_R}{M_1}\right)
      \int_0^{2\pi}h(\theta)H(R,\theta)\,\diff\theta.
\end{aligned}
\end{equation}
Here $[\varepsilon^2]$ denotes the quadratic Taylor coefficient.
To see the cancellation in this formula directly, the
second-order boundary terms in each $M_p(\varepsilon)$
are $X_R^p$ times the corresponding area-variation terms,
since $u_0'(R)=0$, and hence vanish under the fixed-area
constraint.  The constant Fourier mode of the second
harmonic variation contributes $p$ times the same constant
to $\log M_p(\varepsilon)$ and cancels in
$\log M_2-2\log M_1$.

The coefficient in \eqref{eq:neumann-ratio-hessian} is strictly
positive.  First, \eqref{eq:neumann-first-shape-harmonic}
gives
\[
 \int_0^{2\pi}h(\theta)H(R,\theta)\,\diff\theta
   =18\pi R\sum_{m\ne0}\frac{|h_m|^2}{|m|}>0.
\]
Since $X(r)=r(1-r^2)$ is strictly increasing for $0<r<R$,
we have $M_2<X_R M_1$; thus the coefficient of the last
integral in \eqref{eq:neumann-ratio-hessian} is positive.
Moreover, the angular mean of $H(r,\theta)^2$ is
nondecreasing in $r$, by its Fourier expansion.  Weighting
its radial mean by $X^2$ rather than by $X$ therefore
increases its normalized average.  Consequently,
\[
 \frac{\int X^2H^2\,\diff A_{\mathbb H}}{M_2}
 \ge\frac{\int XH^2\,\diff A_{\mathbb H}}{M_1}>0.
\]
Both terms together make the quadratic coefficient positive.
Since $\operatorname{Area}_{\mathbb H}(\Omega_\varepsilon)
=\frac12$ and the first variation vanishes, the same is
true for $\mathcal R_N=1-M_1(\varepsilon)^2/
(\frac12 M_2(\varepsilon))$, as asserted.
\end{proof}

The positive second variation is not uniformly coercive in the
ordinary boundary $L^2$ norm.  To see this explicitly, take
$h_m(\theta)=\cos(m\theta)$ in
Proposition~\ref{prop:neumann-disk-local-minimum}.  Formula
\eqref{eq:neumann-first-shape-harmonic} gives
\[
 H_m(r,\theta)=\frac{9R}{m}\left(\frac rR\right)^m
                  \cos(m\theta),\qquad m\ge1.
\]
Consequently the two bulk terms in
\eqref{eq:neumann-ratio-hessian} are $\Ordo(m^{-3})$,
whereas the boundary term is of order $m^{-1}$.  More precisely,
putting $X_R=R(1-R^2)$ and
$M_p=\int_{|z|<R}\e^{pu_0}\,\diff A_{\mathbb H}$ as above,
we obtain
\begin{equation}\label{eq:neumann-high-mode-hessian}
 c[h_m]=\frac{27}{2m}(1-\rho_*)
   \left(\frac{X_R^2}{M_2}-\frac{X_R}{M_1}\right)
   +\Ordo(m^{-3})
 =\frac{3.7682729\ldots}{m}+\Ordo(m^{-3}).
\end{equation}
Here $M_1=2(\operatorname{arctanh}R-R)$ and $M_2=1/18$.
Thus every nonzero smooth Fourier mode has a positive quadratic
coefficient, but these coefficients tend to zero with frequency.
A proof of local optimality would need an appropriate boundary
norm and control of the nonlinear remainder, not merely pointwise
positivity of the second variation.

\paragraph{A coarea identity and the isoperimetric obstruction.}
There is an exact distribution-function formulation of the
Neumann-cell problem.  Let $E_t=\{z\in\Omega:u(z)<t\}$,
$A(t)=\operatorname{Area}_{\mathbb H}(E_t)$, and let
$\ell(t)$ denote the curvature $-1$ hyperbolic length of the interior
level arcs $\{u=t\}\cap\Omega$.  At a regular level for which the
coarea formula applies, the Neumann condition on the part
of $\partial E_t$ lying in $\partial\Omega$ gives
\[
 \begin{aligned}
 \frac12-A(t)
 &=\frac1{4\pi}\int_{\{u=t\}\cap\Omega}
                   |\nabla_{-1}u|\,\diff s_{-1},\\
 A'(t)&=\frac1{4\pi}\int_{\{u=t\}\cap\Omega}
                   \frac{\diff s_{-1}}{|\nabla_{-1}u|}.
 \end{aligned}
\]
Consequently, Cauchy--Schwarz yields
\begin{equation}\label{eq:neumann-cell-coarea}
 \left(\frac12-A(t)\right)A'(t)
       \ge\frac{\ell(t)^2}{16\pi^2}.
\end{equation}
When $E_t$ lies strictly inside the cell, the ordinary
hyperbolic isoperimetric inequality supplies
$\ell(t)^2\ge (4\pi)^2 A(t)(A(t)+1)$, and hence
\begin{equation}\label{eq:neumann-cell-interior-iso}
 \left(\frac12-A(t)\right)A'(t)
       \ge A(t)(A(t)+1).
\end{equation}
The circular cell has equality.  Importantly, this inequality
alone has the \emph{wrong direction} for the desired lower
bound on squared discrepancy: after using $A$ as the
independent variable, it bounds $\diff t/\diff A$ from
above by the circular value, which by itself tends to make
the distribution of $\e^u$ more uniform.  Thus a naive
isoperimetric rearrangement is not a proof of
Conjecture~\ref{conj:neumann-cell-sharp}.

For the levels which meet the cell boundary, $\ell(t)$ is
only the length of the \emph{interior} level arcs, not the
full perimeter of $E_t$.  The boundary arcs are Neumann
streamlines and may contain saddle--maximum connections.
The local positivity in Proposition~\ref{prop:neumann-disk-local-minimum}
provides evidence for, but does not prove, the global sharp
comparison.  A possible path toward \eqref{eq:neumann-cell-lower} is to
find a sharp two-moment comparison for $\e^u$ which takes
account of the full boundary-contact distribution and the
critical-point topology.  An intermediate step would be to
prove the comparison for a controlled class of finite
Neumann graphs and then extend it by exhaustion.  Only after
such a cell inequality has been proved should one pass to
near-minimizers of the disk discrepancy, controlling multiple
zeros, noncompact cells, and boundary loss.  The numerical
value $\rho_*$ is therefore a conjectural universal lower
bound, not a consequence of \eqref{eq:neumann-cell-coarea}.

A second exact identity may be useful in such a comparison.
For every $p>0$, integration by parts in the Neumann equation,
using the vanishing of $\e^{pu}$ at its logarithmic pole, gives
\begin{equation}\label{eq:neumann-exponential-energy}
 \int_\Omega\e^{pu}|\nabla_{-1}u|^2
   \,\diff A_{\mathbb H}
 =\frac1p\int_\Omega\e^{pu}\,\diff A_{\mathbb H}.
\end{equation}
In particular $X=\e^u$ satisfies the fixed Dirichlet
quotient identity
\[
 \int_\Omega|\nabla_{-1}X|^2\,\diff A_{\mathbb H}
 =\frac12\int_\Omega X^2\,\diff A_{\mathbb H}.
\]
This supplies an additional analytic constraint on the
admissible distribution functions, beyond coarea and the
ordinary isoperimetric inequality.  A sharp argument must
still incorporate the one-zero condition and the geometry
of the Neumann graph.
 Gradient basins of a general differential are
usually irregular Neumann cells whose boundaries connect saddles
and maxima. Their high-level islands may supply the additional
loss compared with the circular model.

The proposed Neumann-cell inequality provides a direct route to
an explicit lower bound for the hyperbolic discrepancy density:
\begin{equation}\label{eq:conjectured-rho-lower}
 \rho_{\mathbb H}\ge\rho_*=0.0522117855\ldots .
\end{equation}
Indeed, on any union of complete gradient-flow basins containing
one simple zero each, the cells have normalized hyperbolic area $\frac12$ and the
conjecture gives this bound cell by cell. Minimizing with one common
scale cannot improve upon minimizing independently on the cells.
To pass from such a union to the asymptotic disk problem defining
$\rho_{\mathbb H}$, one must still justify the gradient-flow
exhaustion for near-minimizers and control exceptional cells and
boundary contributions. The conjecture thus identifies the local
geometric estimate needed for the global lower bound.
The bound would not by itself determine the exact value, since the
periodic construction only gives the upper bound
$\rho_{\mathbb H}\le\rho_{12,3}\simeq0.0523107153$, and the
nonperiodic improvement discussed above indicates that the latter
is not expected to be sharp. Nor would the lower bound alone give
the same numerical improvement for the unrestricted maximal
asymptotic variance in the folding argument.

Other regular finite-area tessellations with the same area per
center are $\{8,4\}$, $\{6,6\}$, $\{5,10\}$ and the ideal-square
$\{4,\infty\}$ tessellation.  Our fundamental-cell computations
for the other compact regular models give
$\rho_{8,4}\simeq0.05316028$,
$\rho_{6,6}\simeq0.05811656$, and
$\rho_{5,10}\simeq0.07676403$.
Thus $\{12,3\}$ has the smallest computed discrepancy among these
four compact regular tessellations.  Cusped quotients require care.
A logarithmic differential with a nonzero residue has infinite
squared hyperbolic discrepancy near a cusp, so the finite-energy
problem is not obtained simply by replacing a compact fundamental
cell with an ideal polygon.

\subsection{Schwarzian variance and zero packing}
\label{subsec:schwarzian-variance-lower}
The discrepancy upper bound from the regular $\{12,3\}$
cell also gives a numerical \emph{lower} bound for maximal
Schwarzian asymptotic variance.  This deduction does not use
the strict finite-scale folding theorem of Section~2, nor
the conjectural Neumann-cell inequality.

For an analytic function $H$ on $\D$, define
\[
 \sigma_2(H)^2:=\limsup_{r\to1^-}
  \frac{1}{L(r)}\int_\D
  |H(rz)|^2(1-|z|^2)^3\,\dA(z),
 \quad L(r)=\log\frac1{1-r^2},
\]
and for a holomorphic function $g$ put
\[
 \sigma(g)^2:=\limsup_{r\to1^-}\frac1{L(r)}
       \int_{\T}|g(r\zeta)|^2\,\diff s(\zeta).
\]
Let $\mathfrak S_{\max}$ denote the corresponding
maximal squared Schwarzian amplitude,
\[
 \mathfrak S_{\max}:=
     \sup_{f\in\mathcal S}\sigma_2(S_f)^2,
 \qquad S_f=(\Nop(f))'-\tfrac12(\Nop(f))^2.
\]

\Needspace{5\baselineskip}
The following proposition relates the maximal Schwarzian asymptotic variance to the hyperbolic zero-packing discrepancy through the infinitesimal Bergman-projection method.

\begin{prop}
\label{prop:schwarzian-lower-from-packing}
With the above normalizations,
\begin{equation}\label{eq:schwarzian-lower-from-packing}
 \mathfrak S_{\max}\ge6(1-\rho_{\mathbb H})
 \ge6(1-\rho_{12,3})>5.686116.
\end{equation}
\end{prop}
\begin{proof}
Let $\mu\in L^\infty(\D)$ with $\|\mu\|_\infty\le1$,
initially vanishing on some disk $|\zeta|<\varepsilon$.
Define an exterior Beltrami coefficient by
\[
 \nu(w)=\begin{cases}
 \bar w^{-2}\mu(1/\bar w),&|w|>1,\\
 0,&|w|\le1.
 \end{cases}
\]
Then $\|\nu\|_\infty\le1$ and $\nu$ has compact support.
For $|\lambda|<1$, let $\Phi_\lambda$ be the principal
solution of $\bar\partial\Phi_\lambda
=\lambda\nu\,\partial\Phi_\lambda$, normalized at
infinity, and put
\[
 f_\lambda(z)=
 \frac{\Phi_\lambda(z)-\Phi_\lambda(0)}
 {\Phi_\lambda'(0)},\qquad z\in\D.
\]
Each $f_\lambda$ belongs to $\mathcal S$ and depends
holomorphically on $\lambda$.  With the ordinary Bergman
projection
\[
 P\mu(z)=\int_\D\frac{\mu(\zeta)}
 {(1-z\bar\zeta)^2}\,\dA(\zeta),
\]
the Cauchy-transform first-variation formula, followed by
the substitution $w=1/\bar\zeta$, yields
\[
 \left.\partial_\lambda\log f_\lambda'(z)
       \right|_{\lambda=0}
   =-P\mu(z)+P\mu(0).
\]
Since $S_{f_0}=0$, it follows that
\begin{equation}\label{eq:schwarzian-beltrami-linearization}
 \left.\partial_\lambda S_{f_\lambda}(z)
       \right|_{\lambda=0}=-(P\mu)''(z).
\end{equation}

Fix $0<k<1$.  Holomorphic dependence on $\lambda$ and
Parseval's inequality on the parameter circle give, at every
$z\in\D$,
\[
 \int_{\T}|S_{f_{k\zeta}}(z)|^2\,\diff s(\zeta)
   \ge k^2 |(P\mu)''(z)|^2.
\]
Integrate with respect to
$L(r)^{-1}(1-|z|^2)^3\,\dA(z)$, evaluating all
functions at $rz$.  To pass to the upper limit in $r$,
use the classical Kraus bound
\[
 (1-|w|^2)^2|S_f(w)|\le6,\qquad f\in\mathcal S.
\]
Indeed, uniformly for $\zeta\in\T$ and $r$ bounded away
from zero, the normalized integral is bounded above by
\[
 \frac{36}{L(r)}\int_0^1
       \frac{(1-t)^3}{(1-r^2t)^4}\,\diff t
 \le\frac{36}{r^2}.
\]
Along a sequence of radii realizing the upper limit of the
right-hand side, reverse Fatou (applied to these uniformly
bounded nonnegative parameter integrals) gives
\[
 \mathfrak S_{\max}\ge k^2\sigma_2((P\mu)'')^2.
\]
For a Bloch function $g$, the elementary asymptotic
differentiation identity, obtained by comparing Taylor
coefficients with the beta integral, reads
\[
 \sigma_2(g'')^2=6\sigma(g)^2.
\]
Since $P\mu$ is a Bloch function, we conclude that
$\mathfrak S_{\max}\ge6k^2\sigma(P\mu)^2$.
Letting $k\to1^-$ and optimizing over $\mu$ gives
\[
 \mathfrak S_{\max}\ge
 6\sup_{\|\mu\|_\infty\le1}\sigma(P\mu)^2
 =6(1-\rho_{\mathbb H}),
\]
by the hyperbolic discrepancy identity
of~\cite{HedenmalmAJM}.  The initial requirement that
$\mu$ vanish near the origin can be removed: its deleted
part has Bergman projection analytic past $\overline\D$
and contributes no asymptotic variance, including in the
cross terms by Cauchy--Schwarz.  Finally apply
Proposition~\ref{prop:12-3-cell-upper-bound}.
\end{proof}

\begin{rem}
Conjecture~\ref{conj:neumann-cell-sharp}, together with
control of gradient-flow exhaustions, would give
$\rho_{\mathbb H}\ge\rho_*$.  The two cell bounds would then
sandwich the \emph{infinitesimal} Schwarzian coefficient:
\[
 5.6861357082\ldots=6(1-\rho_{12,3})
 \le6(1-\rho_{\mathbb H})
 \le6(1-\rho_*)=5.68672928\ldots.
\]
This would not provide an upper bound for the unrestricted
quantity $\mathfrak S_{\max}$, which can contain additional
nonlinear contributions.
\end{rem}

\subsection{An exploratory Beltrami family from the zero lattice}
\label{subsec:candidate-boundary-domains}

The regular $\{12,3\}$ lattice is an especially explicit example:
its optimized cell discrepancy $\rho_{12,3}\simeq0.0523107153$
is very close to the conjectural lower benchmark
$\rho_*=0.0522117855\ldots$.  This agreement is a reason to
investigate its associated Beltrami family, not evidence that a
particular phase maximizes asymptotic variance.  The cell discrepancy
is a quantity attached to the holomorphic one-form, whereas the
variance of its phase projection is a distinct quantity and need not
equal $1-\rho_{12,3}$.

One natural Beltrami coefficient for the exterior conformal problem
is the \emph{conjugate phase} of $f_{12}$,
\[
 \nu_{12}(z)=\begin{cases}
 \overline{f_{12}(z)}/|f_{12}(z)|,&z\in\D,\\
 0,&z\notin\D.
 \end{cases}
\]
The value at a zero is immaterial. For $|\lambda|<1$ let
$\Psi_\lambda$ be the principal solution
\[
 \bar\partial\Psi_\lambda
 =\lambda\nu_{12}\,\partial\Psi_\lambda,
 \qquad \Psi_\lambda(z)=z+\mathrm O(z^{-1})\quad(z\to\infty).
\]
Then $\Omega_\lambda=\Psi_\lambda(\D)$ is an illustrative image
domain and $\Psi_\lambda|_{\D_{\mathrm e}}$ is its exterior
conformal chart.  In contrast to the Teichm\"uller coefficient
$\overline{f_{12}}/f_{12}$ attached to a quadratic differential,
$\nu_{12}$ is the conjugate phase corresponding to the Bergman
duality test below.  No optimality assertion for this fixed Beltrami
coefficient or for its nonlinear deformations is intended.

To make the connection with the Bergman duality of
\cite{HedenmalmAJM} explicit, let
\[
 \phi_{12}(z):=\frac{f_{12}(z)}{|f_{12}(z)|}
   =\overline{\nu_{12}(z)},\qquad
 g_{12}:=\Bergman[\phi_{12}],
 \quad z\in\D,
\]
where the ordinary Bergman projection is
\[
 \Bergman[h](z):=\int_\D
       \frac{h(\zeta)}{(1-z\bar\zeta)^2}\,\dA(\zeta).
\]
Thus the phase in the Bergman projection is the conjugate of the
Beltrami coefficient. For the principal solution above,
Cauchy-transform differentiation at $\lambda=0$ gives, for $|w|>1$,
\begin{equation}\label{eq:12-3-beurling-bergman-reflection}
\begin{aligned}
 \left.\partial_\lambda\log\Psi_\lambda'(w)
                          \right|_{\lambda=0}
 &=-\int_\D\frac{\nu_{12}(z)}{(w-z)^2}\,\dA(z)\\
 &=-\frac1{w^2}\,
      \overline{g_{12}\!\left(\frac1{\bar w}\right)}.
\end{aligned}
\end{equation}
Consequently, for $0<r<1$,
\begin{equation}\label{eq:12-3-phase-variance-identity}
 \int_\T\left|
  \left.\partial_\lambda\log
   \Psi_\lambda'\!\left(\frac\zeta r\right)
  \right|_{\lambda=0}\right|^2\,\diff s(\zeta)
  =r^4\int_\T|g_{12}(r\zeta)|^2\,\diff s(\zeta).
\end{equation}
Since $r^4\to1$, the infinitesimal exterior logarithmic
asymptotic variance of this Beltrami family equals
$\sigma(g_{12})^2$.  This exact identification does not determine
that variance numerically, and it does not identify it with the
cell value $1-\rho_{12,3}$.  Replacing $\nu_{12}$ by its conjugate
would instead produce $\Bergman[\nu_{12}]$ in
\eqref{eq:12-3-beurling-bergman-reflection} and, in general,
a different variation.

Writing $\mathcal C$ for the Cauchy transform,
$\mathcal B=\partial\mathcal C$ for the Beurling transform, and
$M_\nu$ for multiplication by $\nu_{12}$, the principal solution
has the norm-convergent Neumann expansion in the Beltrami density
for $|\kappa|<1$:
\begin{equation}\label{eq:beltrami-neumann-candidate}
 \Psi_\kappa(z)=z+
 \sum_{n\ge1}\kappa^n
 \mathcal C(M_\nu\mathcal B)^{n-1}\nu_{12}(z).
\end{equation}
The first-order exterior logarithmic derivative is
$\mathcal B\nu_{12}$, equivalently the reflected conjugate of
$g_{12}$ by \eqref{eq:12-3-beurling-bergman-reflection}.
This makes the family a concrete test of the AJM duality mechanism,
without implying infinitesimal or finite-$\kappa$ optimality.
Preliminary numerical circle-mean experiments gave values broadly
consistent with a variance near $0.95$, but logarithmic convergence
is slow, and no exact limiting value is claimed.

\begin{figure}[t]
 \centering
 \includegraphics[width=.97\linewidth]{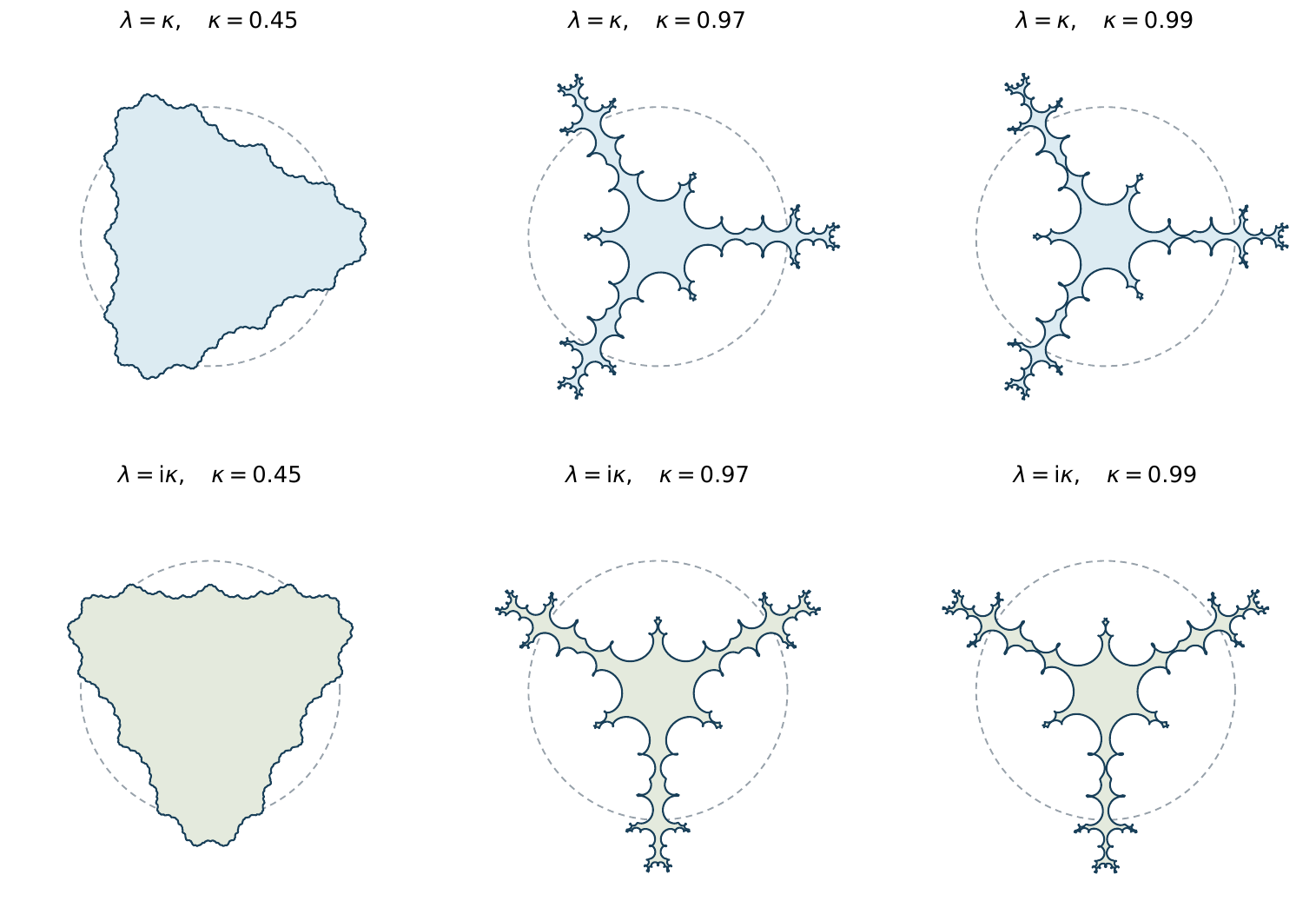}
 \caption{Candidate domains for $\kappa=0.45$, $0.97$, and $0.99$
 (columns, from left to right). The top row corresponds to
 $\lambda=\kappa$ and the bottom row to
 $\lambda=\imag\kappa$, displaying the Y-shaped orientation.
 The imaginary-phase curves are obtained by the exact rotations
 in \eqref{eq:phase-rotated-domains} rather than separate Beltrami solves.
 The dashed circle is the original unit circle. The contours are
 finite-resolution numerical approximations.}
 \label{fig:candidate-domain-phases}
\end{figure}

\paragraph{Exact rotation covariance.}
Write $\zeta=\exp(\pi\imag/6)$. Since the form $f_{12}$ has a
simple zero at the order-twelve elliptic center, its normalized local
expression satisfies
\[
 f_{12}(\zeta z)=\zeta f_{12}(z),\qquad
 \nu_{12}(\zeta z)=\zeta^{-1}\nu_{12}(z).
\]
Accordingly $g_{12}(\zeta z)=\zeta g_{12}(z)$, so the Taylor
coefficients of $g_{12}$ vanish outside the congruence class
$n\equiv1\pmod{12}$. This is consistent with the choice of
conjugate phase in the Beltrami equation.
Apply the coordinate change $z\mapsto\zeta z$ to the principal
Beltrami solution. The transformation rule for its complex dilatation,
together with the principal normalization at infinity, gives the
\emph{exact} identity
\begin{equation}\label{eq:phase-rotation-covariance}
 \Psi_{\lambda\zeta^{-3}}(z)
 =\zeta^{-1}\Psi_\lambda(\zeta z),\qquad |\lambda|<1.
\end{equation}
In particular, for $0<\kappa<1$,
\begin{equation}\label{eq:phase-rotated-domains}
 \Omega_{-\kappa}=\exp(-\pi\imag/3)\Omega_\kappa,
 \qquad
 \Omega_{\imag\kappa}=\exp(\pi\imag/6)\Omega_\kappa.
\end{equation}
Since $\zeta^4=\exp(2\pi\imag/3)$ and $\zeta^{-12}=1$,
all $\Omega_\lambda$ are invariant under rotation through
$2\pi/3$. Thus the choices $\lambda=-\kappa$ and
$\lambda=\imag\kappa$ produce no new shapes up to rotation.
Figure~\ref{fig:candidate-domain-phases} compares the real and imaginary
phases, which differ by an exact rotation. The negative real phase,
being another rotation, is omitted from the figure.

\paragraph{Numerical interpretation near $\kappa=1$.}
Figure~\ref{fig:candidate-domain-phases} shows increasingly fine
boundary structure as $\kappa$ approaches $1$. The computations use
polar grids of $1152\times4608$ points for $\kappa=0.45$ and
$384\times1536$ points for $\kappa=0.97,0.99$.
The contours are qualitative: a small fixed-point residual does not
certify the finest corrugations, a limiting fractal dimension, or
spectral optimality. The threefold symmetry is exact.

\subsection{Toward limiting boundary geometry}
\label{subsec:change-with-k}

The small-dilatation expansion is universal:
\begin{equation}\label{eq:small-k-variance-selection}
 \Sigma^2(\kappa)
 =\bigl(1-\rho_{\mathbb H}\bigr)\kappa^2
 +o(\kappa^2),\qquad \kappa\to0^+.
\end{equation}
It does not select an extremizer at finite $\kappa$.
As $\kappa$ grows, the candidate may remain stable, bifurcate,
or be superseded by nonperiodic and degenerating configurations.
At $\kappa=1$ the Beltrami bound no longer guarantees a
quasiconformal homeomorphism, and continuity of
$\Sigma^2(\kappa)$ at that endpoint is unknown.

\paragraph{Exploratory limiting geometry.}
Random conformal snowflakes~\cite{BeliaevSmirnovSnowflakes}
and Makarov's multifractal theory~\cite{MakarovFineStructure}
motivate examining possible boundary limits, without identifying
the fixed $\{12,3\}$ phase as an optimizer. Write
$\D_{\mathrm e}=\{|z|>1\}$ and
$\psi_\kappa=\Psi_\kappa|_{\D_{\mathrm e}}$.
Normal-family compactness gives subsequential locally uniform
limits of $\psi_\kappa$ as $\kappa\to1^-$; uniqueness and
spectral optimality are unknown.

Let $\Sigma$ be the class of normalized exterior maps
$\psi(z)=z+\mathrm O(z^{-1})$. Define
\[
 \beta_\psi(t):=\limsup_{R\to1^+}
 \frac{\log\displaystyle\int_{\T}
    |\psi'(R\zeta)|^t\,\diff s(\zeta)}
 {\log\frac1{R-1}},\qquad
 B_\Sigma(t):=\sup_{\psi\in\Sigma}\beta_\psi(t).
\]
The issue is to determine the limiting spectra and whether
perturbations of the zero configuration improve them.
The figures alone imply neither uniqueness nor fractality.

\paragraph{A possible driver--environment mechanism.}
The exact Littlewood--Paley identity for exponentials of the
logarithmic derivative suggests that the squared nonlinearity acts
as a driver, while the exponential of a harmonic function supplies
an exponent-dependent environment. This leads to a weighted form
of hyperbolic zero packing. Its possible invariance under suitable
log-harmonic perturbations may help explain why a single lattice
could be relevant for several exponents. At present this is a
heuristic interpretation, not an additional ingredient in the
strict variance proof or a proof of spectral optimality.

\subsection{Questions suggested by the figures}
\label{subsec:questions-extremal-domains}

Open questions include the Neumann-cell comparison
\eqref{eq:neumann-cell-lower}, control of boundary-touching
gradient cells, solutions of the limiting Euler--Lagrange equation
\eqref{eq:limiting-zero-packing-fixed-point}, and stability of
periodic zero configurations on larger covers. The geometry of
possible limits near $\kappa=1$ presents a separate compactness
problem. Stationarity on a fixed quotient does not imply universal
optimality.

\medskip
\noindent\emph{Numerical details.} An a posteriori enclosure of the
$\{12,3\}$ cell ratio is supplied as a numerical supplement,
including bounds for the Neumann residual and moment quadrature. It yields
$\rho_{12,3}<0.052314$ and hence $6(1-\rho_{\mathbb H})>5.686116$.

\medskip

\noindent\textbf{Acknowledgements.}
The author thanks Oleg Ivrii for sharing his ideas and unpublished work.
The author acknowledges support from the Ministry of Science and Higher
Education of the Russian Federation under agreement 075-15-2025-013.
The author acknowledges the assistance of OpenAI's ChatGPT in
expository preparation, exploratory computations, and preliminary checks.
The author takes responsibility for the mathematical statements and their
verification.

\end{document}